\documentclass[10pt, a4paper]{scrartcl}
\usepackage[T1]{fontenc} 
\usepackage[english]{babel}

\usepackage{amsmath}
\usepackage{amssymb}
\usepackage{amsthm}

\newcounter{count_lemma}
\newcounter{count_theorem}

\newtheorem{lemma}[count_lemma]{Lemma}
\newtheorem{theorem}[count_theorem]{Theorem}

\newcommand{\R}{\mathbb{R}}
\newcommand{\C}{\mathbb{C}}

\newcommand{\N}{\mathbb{N}}
      
\newcommand{\dollar}{\hfill \fbox{}\\}
\newcommand{\transpose}{{\hspace{-0.01cm}\scriptscriptstyle\mathrm T}} 
\newcommand{\Oc}{\mathcal O}
\newcommand{\Kc}{\mathcal K}

\newcommand{\Rc}{\mathcal R}

\newcommand{\Wc}{\mathcal W}
\newcommand{\Nc}{\mathcal N}
\newcommand{\Nct}{\tilde{\mathcal N}}

\newcommand{\bfj}{\textbf{j}}
\newcommand{\bfk}{\textbf{k}}

\newcommand{\bfe}{\textbf{e}}
\newcommand{\bfm}{\textbf{m}}
\newcommand{\bfr}{\textbf{r}}

\newcommand{\T}{\mathfrak{T}}

\newcommand{\D}{\mathfrak{D}}
\usepackage{tikz}
    \usetikzlibrary{arrows}
    \usetikzlibrary{trees}     
\usepackage{floatrow}
\usepackage[margin=0cm,normal,bf,small]{caption}%
\usepackage{wrapfig} 
\usepackage{tabu}
\usepackage{multirow}
\usepackage{multicol}
\usepackage{verbdef}
\usepackage[vlined,ruled]{algorithm2e}
 \newcommand\BlockIf[1]{#1}
 \newcommand\BlockElseIf[1]{#1}
 \newcommand\BlockElse[1]{#1}
\makeatletter
\renewcommand{\uIf}[2]{\If@ifthen{#1}\If@noend{\BlockIf{#2}}}
\renewcommand{\uElseIf}[2]{\ElseIf@elseif{#1}\If@noend{\BlockElseIf{#2}}}
\renewcommand{\Else}[1]{\Else@else\If@endif{\BlockElse{#1}}}
\makeatother

\usepackage{xcolor}
\definecolor{color1}{rgb}{.8,.8,.8}
\definecolor{color2}{rgb}{.65,.65,.65}
\definecolor{color3}{rgb}{.35,.35,.35}
\definecolor{color4}{rgb}{0,0,0}

\def\mh{{\text{e}\hbox{-}}}
\def\mhp{{\text{e}\hbox{+}}}
\begin{document}
\begin{flushright}
March 4, 2015
\end{flushright}

\begin{center}
\textbf{\large{A fast matrix-free algorithm for spectral approximations to the Schr\"odinger equation}}
\end{center}
\begin{center}
\textbf{ \large {Bernd Brumm}} \\
Mathematisches Institut, Universit\"at T\"ubingen, Auf der Morgenstelle 10, D-72076 Germany.\footnote {E-mail: \texttt{brumm@na.uni-tuebingen.de}.}
\end{center}
\vskip 0.5cm 
\textbf{Abstract.} We consider the linear time-dependent Schr\"odinger equation with a time-dependent smooth potential on an unbounded domain. 
A Galerkin spectral method with a tensor-product 
Hermite basis is used as a discretization in space. Discretizing the resulting ODE for the Hermite expansion coefficients 
involves the computation of the action of the Galerkin matrix on a vector in each time step. 
We propose a fast algorithm for the direct computation of this matrix-vector product without actually 
assembling the matrix itself. The costs scale linearly in the size of the basis.   
Together with the application of a hyperbolically reduced basis, this reduces the 
computational effort considerably and helps cope with the infamous curse of dimensionality. 
The application of the fast algorithm is limited to the case 
of the potential being significantly smoother than the solution. 
The error analysis is based on a binary tree 
representation of the three-term recurrence relation 
for the one-dimensional Hermite functions. 
The fast algorithm constitutes an efficient tool for schemes involving the action of a 
matrix due to spectral discretization on a vector, and it is also applicable in the context of spectral approximations for linear problems 
other than the Schr\"odinger equation.

%%%%%%%%%%%%%%%%%%%%%%%%%%%%%%%%%%%%%%%%%%%%%%%%%%%%%%%%%%%%%%%%%%%%%%%%%%%%%%%%%%%%%%
%%%%%%%%%%%%%%%%%%%%%%%%%%%%%%%%%%%%%%%%%%%%%%%%%%%%%%%%%%%%%%%%%%%%%%%%%%%%%%%%%%%%%%
\section*{Introduction}\label{section:Introduction}
We consider the linear time-dependent Schr\"odinger equation
\begin{align}\label{align:SchroedingerEquation}
 i \frac{\partial}{\partial t}\psi(x,t)\!=\!(H\psi)(x,t)\!=\!-\frac{1}{2}(\Delta \psi)(x,t)\!+\!V(x,t)\psi(x,t),\,\, x\!=\!(x_1,\dots,x_N)\!\in\!\R^N
\end{align}
in $N$ spatial dimensions with a possibly time-dependent multiplicative potential $V$ that meets certain regularity conditions 
on a cube $\Omega = [-L,L]^N$ and a solution $\psi(\cdot,t)$ that is essentially supported within $\Omega$, for all times $t \in [0,T]$. For an underlying geometry as simple as in \eqref{align:SchroedingerEquation}, 
spectral methods are a natural means of discretization in space. In a naive approach, the resulting ODE system grows exponentially in $N$, making an accurate approximation 
practically unfeasible even for moderate choices of $N$. For this difficulty, the catch phrase \emph{curse of dimensionality} has been coined. 
Time propagation typically requires computing the action of the Galerkin matrix on a vector in each step, and, in case of a time-dependent potential, the matrix has to be re-assembled.\\%\vspace{\baselineskip}\\
%%%
A promising strategy is a suitable reduction of the spectral approximation basis. E.g., \cite{G07a, G07b} study a spectral approach with collocation on a sparse grid 
in case of a time-\emph{in}dependent potential and periodic boundary conditions with a hyperbolically reduced tensor-product Fourier 
basis. As is pointed out in \cite{L08}, Chapter III.1.4, unlike on a full grid, the resulting coefficient ODE does not exhibit 
a Hermitian matrix, which possibly gives rise to numerical troubles and limits the range of applicable time-stepping methods. As a remedy,  %discrete Fourier Galerkin methods with a simplified mass matrix or a redefined sparse grid are proposed as well as 
a Fourier Galerkin method in combination with an approximation of the potential by a trigonometric polynomial is proposed. We adopt this basic idea from the simpler setting of a periodic problem.\\%\vspace{\baselineskip}\\
%%%
In the present paper, we also employ a spectral Galerkin approach with a reduced basis in combination with a polynomial approximation of (parts of) the potential, but we consider an unbounded domain instead of a periodic problem. Hermite 
functions are a natural and, thus, widely-used spectral basis for the Schr\"odinger equation on unbounded domains, see, e.g., \cite{L08}, Chapter III.1, \cite{FG09} for the linear and \cite{G11} 
for a nonlinear case. Furthermore, we allow the potential to be time-dependent.\\%\vspace{\baselineskip}\\
%%%
Besides basis reduction and potential approximation, we develop a fast algorithm for the direct (i.e., matrix-free)  computation of the aforementioned matrix-vector product 
that further speeds up propagation in time considerably. The basic idea for the fast algorithm was proposed in  \cite{FGL09} in the context of a splitting procedure for the linear Schr\"odinger equation in the semi-classical regime:
One uses a recurrence relation for the univariate Hermite functions and orthogonality to define (never actually assembled) 
coordinate matrices for each coordinate direction that act directly on vectors; these matrices are then formally inserted into the polynomially 
approximated potential. In the present paper, we start from the fact that 
this is equivalent to a suitable entrywise approximation of the Galerkin matrix by Gauss--Hermite quadrature, 
which is briefly derived as is commonly done in the context of Discrete Variable Representations, see \cite{LC00}. 
However, this is only true if the matrices are indexed over a full set of multi-indices. 
We give a detailed analysis for the resulting quadrature error as well as for the error due to a hyperbolical index reduction based on binary tree representations. Both errors are well-behaved if the potential can be sufficiently well approximated 
by a multivariate polynomial. If so, we get bounds $C(\Rc,W,L)K^{-\beta}$ and $C(N,\Rc,W,\beta,L)K^{-\beta}$ for the errors due to quadrature and grid reduction, respectively. Here, $W$ is the part of the potential $V$ that is approximated over an $N$-dimensional index set $\Rc(R)$ with maximal univariate polynomial degree $R$, 
$K$ is the maximal number of basis functions employed in each coordinate direction in the Galerkin approximation, and the coefficients of the approximate solution 
exhibit a decay of order $\beta$ with increasing index.\\%\vspace{\baselineskip}\\
%%%
The fast algorithm has the following advantages: First of all, it scales only linearly in the size of the basis. 
In addition, it allows for more general kinds of index reductions than the hyperbolically reduced index set considered in the present work. 
Our approach avoids quadrature at all. In contrast, in the chemical literature, 
there is a matrix-free approach based on ingenious sequential summations for the matrix-vector product that employs a (still exponentially large) reduced basis 
and treats the problem of huge quadrature-grids using a 
nonproduct Smolyak grid quadrature, see \cite{AC09, AC11a, AC11b, AC12}. 
Finally, our algorithm constitutes a useful tool for a variety of time integration schemes involving Galerkin matrix-vector products, and 
the strategy is generic and, thus,  
serviceable also for spectral Galerkin approximations to other linear problems using a Galerkin basis of 
algebraic orthogonal polynomials. 
An article on the treatment of second-order partial differential equations on bounded domains 
with different kinds of boundary conditions involving a Legendre Galerkin basis is currently in preparation.\\
%%%
The limitations of the fast algorithm are as follows: The method is applicable only if the potential is significantly smoother than the solution, i.e., $K \gg R$, wherever the solution does not essentially vanish.
In the present work, we use it in the context of a Lanczos-based time propagation scheme, which yields an error due to a perturbation of the Lanczos process. 
This error may become dominant in cases where the size of the basis is not sufficiently large. 
Finally, our presentation is limited to the setting where the solution stays within a cube and basis functions localized around zero are applicable. This restriction is 
dispensable, see \cite{FGL09} for an application of the fast algorithm with a moving wavepacket basis. An adaptation of the error analysis as given in the present work to the case of an evolving basis is possible, but complicated and technically more involved. 
For the sake of readability, we therefore restrict the scope of the present work to this somewhat idealized setting.

%%%
In Section \ref{section:SemiDiscretizationInSpace}, we deduce 
the ODE system for the Hermite expansion coefficients from the Galerkin ansatz with a reduced index set and a polynomially approximated potential. 
Section \ref{section:SemiDiscretizationInTime} briefly outlines the discretization in time by Magnus integrators, where the matrix exponential is approximated using the Lanczos method. 
Section \ref{section:FastAlgorithm} contains the fast algorithm for the matrix-free computation of the action of the Galerkin matrix on a vector in each 
Lanczos step and illustrates the computational speed-up. 
The connection between the algorithm and Gauss--Hermite quadrature for the Galerkin matrix is shown in Section \ref{section:GaussHermiteQuadrature}. 
A detailed error analysis is given in Section \ref{section:ErrorAnalysis}. 
Section \ref{section:NumericalExperiments} presents some numerical experiments 
confirming the theoretical results.
%%%%%%%%%%%%%%%%%%%%%%%%%%%%%%%%%%%%%%%%%%%%%%%%%%%%%%%%%%%%%%%%%%%%%%%%%%%%%%%%%%%%%%
%%%%%%%%%%%%%%%%%%%%%%%%%%%%%%%%%%%%%%%%%%%%%%%%%%%%%%%%%%%%%%%%%%%%%%%%%%%%%%%%%%%%%%
\section{Semi-discretization in space}\label{section:SemiDiscretizationInSpace}
%--------------------------
\subsection{Hermite basis}\label{subsection:ConstructionHermite}
Starting from $\varphi_{-1} \equiv 0$ and $\varphi_0(x) = \pi^{-1/4} e^{-x^2/2}$, the three-term reccurrence relation 
\begin{align}\label{align:Recurrence1D}
 x\varphi_k(x) = \sqrt{\frac{k+1}{2}}\varphi_{k+1}(x) + \sqrt{\frac{k}{2}}\varphi_{k-1}(x),\qquad k \geq 0,
\end{align}
yields a complete $L^2(\R)$-orthonormal set $\{\varphi_k\}_{k \in \N}$ of Schwartz functions, in particular, $(\varphi_j,\varphi_k) = \delta_{jk}$, 
where $(f,g) = \int f \overline{g}$ denotes the standard $L^2$-inner product. An explicit expression is $\varphi_k(x) = \pi^{-1/4}\left(2^k k!\right)^{-1/2}H_k(x)e^{-x^2/2}$,
where $H_k$ denotes the classical Hermite polynomial of degree $k$. The Hermite functions are readily seen to be the 
eigenfunctions of the harmonic oscillator, i.e.,
\begin{align}\label{align:HermiteEigenfunctions1D}
 \frac{1}{2}(p^2 + q^2)\varphi_k = \left(k + \frac{1}{2}\right)\varphi_k,
\end{align}
where $(q \psi)(x) = x \psi(x)$ and $(p \psi)(x) = -i \frac{\partial}{\partial x}\psi(x)$ denote the position and momentum operators, respectively, 
see, e.g., \cite{AS65}, Section 22, \cite{L08}, Chapter III.1.1, or \cite{Th00}, Section 7.7, for the construction of the Hermite basis and its properties.
In higher dimensions, we consider tensor-products of Hermite functions, i.e.,
\begin{align*}
 \varphi_\bfk(x) = \varphi_{k_1}(x_1)\dots \varphi_{k_N}(x_N),
\end{align*}
where $\bfk = (k_1,\dots,k_N) \in \N^N$ is a multi-index and $\varphi_{k_l}$ are univariate Hermite functions as above, $1 \leq l \leq N$. 
Again, $\{\varphi_\bfk\}_{\bfk \in \N^N}$ is a complete $L^2(\R^N)$-orthonormal set of functions. Due to the eigenfunction property \eqref{align:HermiteEigenfunctions1D}, we find
\begin{align}\label{align:EigenfunctionRelation}
 \frac{1}{2}\sum_{l=1}^N \left(p_l^2 + q_l^2\right)\varphi_\bfk = \frac{1}{2}\left(-\Delta + \sum_{l=1}^N q_l^2\right)\varphi_\bfk = \sum_{l=1}^N \left(k_l + \frac{1}{2}\right)\varphi_\bfk,
\end{align}
where $(q_l\psi)(x) = x_l \psi(x)$ and $(p_l\psi)(x) = -i \frac{\partial}{\partial x_l}\psi(x)$.
%--------------------------
 \begin{figure}
 \floatbox[{\capbeside\thisfloatsetup{capbesideposition={right,top},capbesidewidth=4cm}}]{figure}[\FBwidth]
 {\caption{Univariate Hermite functions for some choices of $k$. $\varphi_k$ is even if $k$ is even, otherwise odd. The largest extremal is bounded by $\sqrt{2(k+1)}$.}
  \label{figure:HermiteFunctions}}
 {\includegraphics[trim = .7cm .2cm 1.4cm .3cm,clip,scale=0.65]{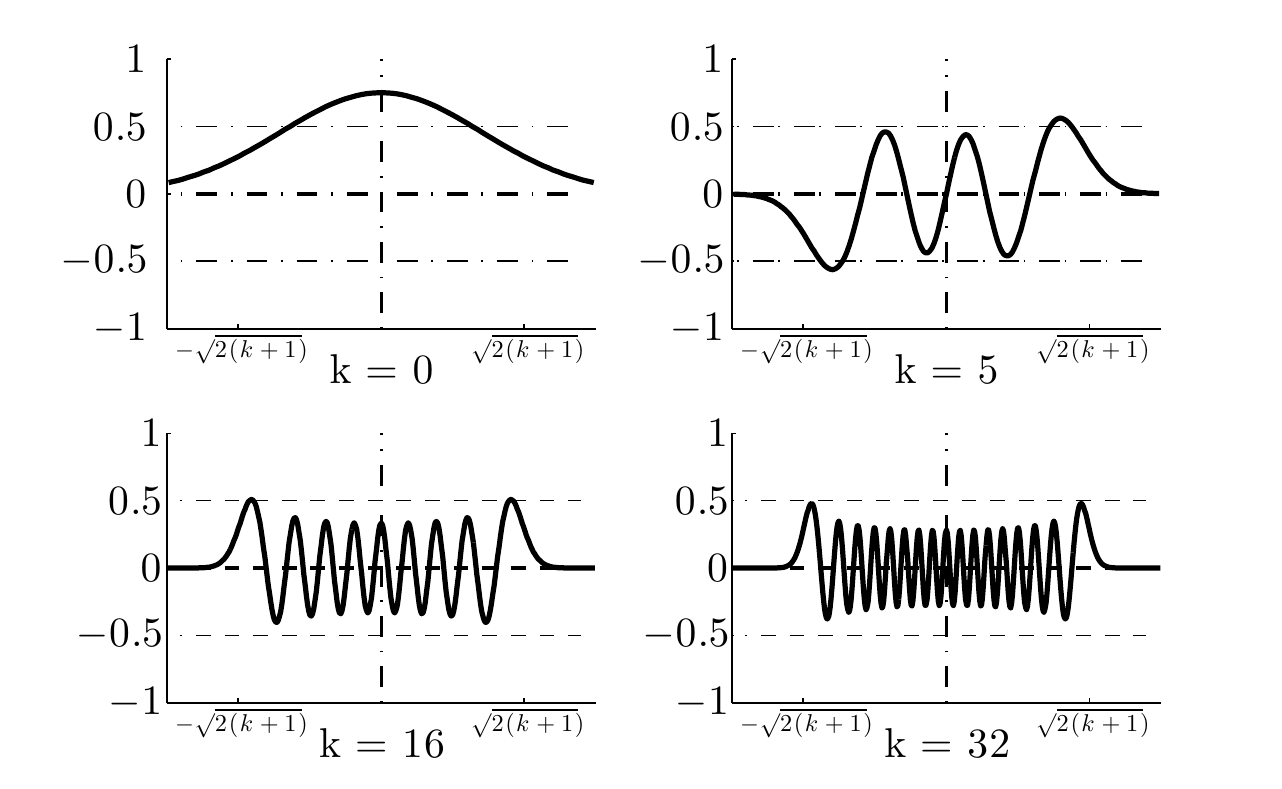}}
\end{figure}\\
%--------------------------
%%%%%%%%%%%%%%%%%%%%%%%%%%%%%%%%%%%%%%%%%%%%%%%%%%%%%%%%%%%%%%%%%%%%%%%%%%%%%%%%%%%%%%
\subsection{Galerkin ansatz with reduced index set}\label{subsection:GalerkinAnsatz}
An approximation
\begin{align}\label{align:ApproximateSolutionLinComb}
\begin{aligned}
 \psi_\Kc(x,t) &= \sum_{\bfk \in \Kc} c_\bfk(t)\varphi_\bfk(x) \in \text{span}\left\{\varphi_\bfk\,\left|\, \bfk \in \Kc\right.\right\} \subseteq L^2(\R^N)
 \end{aligned}
\end{align}
on a finite-dimensional subspace is determined such that
\begin{align}\label{align:GalerkinCondition}
\begin{aligned}
 \left(i \frac{\partial}{\partial t}\psi_{\Kc} - H\psi_{\Kc}, \varphi_\bfj\right) = 0,\qquad \forall\, \bfj \in \Kc,
 \end{aligned}
\end{align}
where $\Kc$ is a multi-dimensional index set
\begin{align}\label{align:FullIndexSet}
 \Kc \subseteq \Kc_{\text{full}} = \left\{\left.\bfk = (k_1,\dots,k_N) \in \N^N\,\right|\, 0 \leq k_l \leq K\right\}
\end{align}
 with at most $K+1$ indices in each direction. If necessary, we write $\Kc = \Kc(K)$ to emphasize the bound $K$ for components $k_l$ of any $\bfk \in \Kc$. 
 Abbreviating $c(t) = (c_\bfk(t))_{\bfk \in \Kc}$ and inserting the ansatz \eqref{align:ApproximateSolutionLinComb} into \eqref{align:GalerkinCondition} yields 
a linear system of ODEs
\begin{align*}
 i \dot{c}(t) = \mathcal{H}_\Kc(t) c(t).
\end{align*}
Furthermore, the eigenfunction relation \eqref{align:EigenfunctionRelation} yields a decomposition ($\bfj, \bfk \in \Kc$)
\begin{align*}
 (\mathcal{H}_\Kc)_{\bfj\bfk} &= \left(\varphi_\bfj, H\varphi_\bfk\right) = \left(\varphi_\bfj, \frac{1}{2}\left(-\Delta + \sum_{l=1}^N q_l^2\right)\varphi_\bfk\right) + \left(\varphi_\bfj, \left(V - \frac{1}{2}\sum_{l=1}^N q_l^2\right)\varphi_\bfk\right)\\
  &= \sum_{l=1}^N \left(k_l + \frac{1}{2}\right) \delta_{\bfj\bfk} + \left(\varphi_\bfj, W\varphi_\bfk\right) = (\mathcal{D}_{\Kc})_{\bfj\bfk} + (\mathcal{W}_\Kc)_{\bfj\bfk},
\end{align*}
where $\mathcal{D}_\Kc = \text{diag}_{\bfk \in \Kc}\left(\sum_{l=1}^N \left(k_l + \frac{1}{2}\right)\right)$ is a diagonal matrix and $(\mathcal{W}_\Kc)_{\bfj\bfk} = (\varphi_\bfj,W\varphi_\bfk)$ stems from a multiplicative potential $W(x,t) = V(x,t) - \frac{1}{2}\sum_{l=1}^N x_l^2$.  
Hence,
\begin{align}\label{align:ODECoefficients}
  i \dot{c}(t) = \mathcal{D}_\Kc c(t) + \mathcal{W}_\Kc(t) c(t).
\end{align}
%--------------------------
\begin{wrapfigure}{r}{0.5\textwidth}
\vspace{-.6cm}
%\hspace{-2cm}
\begin{minipage}{.5\textwidth}
 \includegraphics[trim = .1cm 0cm .1cm 0cm,clip,scale=0.6]{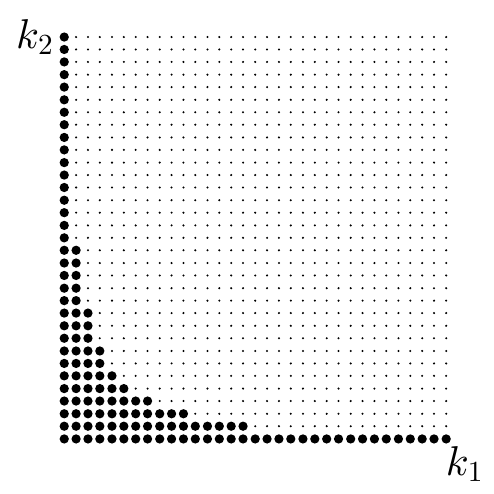}
\end{minipage}%\hspace{1.5cm}
\begin{minipage}{.5\textwidth}
 \includegraphics[trim = .1cm 0cm 0cm 0cm,clip,scale=0.65]{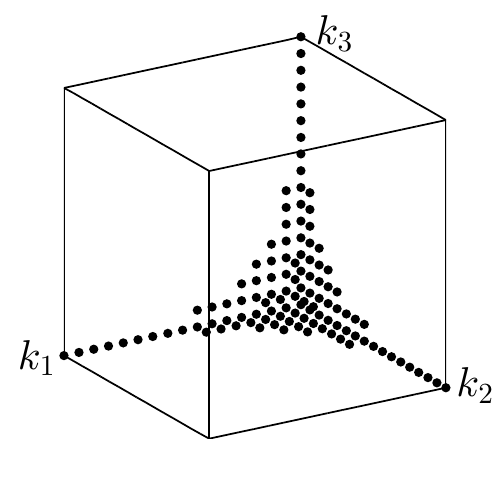}
\end{minipage}
\caption{Hyperbolically reduced index set. Left: $N\!\!=\!\!2$, $K\!\!=\!\!32$. Right: $N\!\!=\!\!3$, $K\!\!=\!\!16$.} 
\label{figure:HyperbolicCross}
\vspace{-.45cm}
\end{wrapfigure}
%--------------------------
In case $\Kc = \Kc_{\text{full}}$, the system \eqref{align:ODECoefficients} consists of $|\Kc| = (K+1)^N$ equations. For growing $N$ and $K$ being only moderate, 
this is not feasible for time integration that requires assembling the matrices $\mathcal{D}_\Kc$ (once) and $\mathcal{W}_{\Kc}(t)$ (in each step) and multiplying them with a vector. Thus, the index set needs to be reduced. We study a hyperbolically
reduced index set
\begin{align*}
 \Kc = \left\{\bfk = (k_1,\dots,k_N)\,\left|\, k_l \geq 0, \prod_{l=1}^N\right. (1 + k_l) \leq K + 1 \right\},
\end{align*}
see the illustration in Figure \ref{figure:HyperbolicCross}. The number of indices employed shrinks to $\left|\Kc\right| = \Oc(K \ln(K)^{N-1})$, see \cite{BG04}. 
Approximating a sufficiently regular function by a Hermite tensor-product expansion over a hyperbolically reduced index set still 
gives a decent approximation under certain regularity assumptions, see \cite{L08}, Thm. III.1.5. More generally, the proposed method allows for any index set $\Kc$ such that, 
if $\bfk \in \Kc$ with $k_l > 0$, then $\bfk - \bfe_l \in \Kc$, 
for all $1 \leq l \leq N$, where $\bfe_l$ is the $l$th unit vector.
%%%%%%%%%%%%%%%%%%%%%%%%%%%%%%%%%%%%%%%%%%%%%%%%%%%%%%%%%%%%%%%%%%%%%%%%%%%%%%%%%%%%%%
\subsection{Approximation of the potential}\label{subsection:ApproximationPotential}
A basic assumption for our fast algorithm is that the potential $W$ can be sufficiently well approximated by an 
interpolation polynomial $W^{\text{pol}}$ 
on a given cube $\Omega = [-L,L]^N$ indexed over a set $\Rc(R) \subseteq \N^N$ with only few nodes, i.e., $|\Rc| \ll |\Kc|$, hence, $R \ll K$.
See Section \ref{subsection:AssumptionOnSolution} for the definition of $L$. 
We consider Chebyshev interpolation, i.e.,
\begin{align*}
 W(x,t) \approx W^{\text{pol}}(x,t) = \sum_{\bfr \in \Rc} \alpha_\bfr(t) T_\bfr(x/L) =  \sum_{\bfr \in \Rc} \alpha_\bfr(t) \prod_{l=1}^N T_{r_l}(x_l/L)
\end{align*}
with coefficients $\alpha_{\bfr}(t)$ that depend on $L$, but not on $K$. The univariate functions $T_{r_l}$ are the Chebyshev polynomials of the first kind that obey the recurrence relation
\begin{align}\label{align:ChebyshevRecurrence}
 \begin{split}
  T_0(x) &= 1,\quad T_1(x) = x,\\
  T_{k+1}(x) &= 2xT_k(x) - T_{k-1}(x),\quad k \geq 1,
 \end{split}
 \begin{split}
  x \in [-1,1].
 \end{split}
\end{align}
In case of a full product grid $\Rc = \Rc_{\text{full}}$, we have exponentially 
decaying coefficients $\alpha_{\bfr}$, see, e.g., \cite{CQHZ06} for a detailed theory of approximation by orthogonal polynomials. 
In case of a reduced index set $\Rc$, a fast polynomial transform on sparse grids can 
be employed to compute the interpolation coefficients efficiently, see, e.g., \cite{CJX14} for an $\Oc(R \ln(R)^{N+1})$ algorithm with still 
exponentially decaying coefficients $\alpha_{\bfr}(t)$ (for $W$ being sufficiently regular) with a hyperbolically reduced set $\Rc$. 
Due to $R \ll K$, the additional costs of doing polynomial interpolation in each time step are negligible. 
In place of \eqref{align:ODECoefficients}, this yields a coefficient ODE
\begin{align}\label{align:ODECoefficientsPol}
   i \dot{c}_{\text{pol}}(t) = \mathcal{D}_{\Kc} c_{\text{pol}}(t) + \mathcal{W}_{\Kc,\text{pol}}(t) c_{\text{pol}}(t).
\end{align}
%%%%%%%%%%%%%%%%%%%%%%%%%%%%%%%%%%%%%%%%%%%%%%%%%%%%%%%%%%%%%%%%%%%%%%%%%%%%%%%%%%%%%%
\subsection{Assumption on the solution}\label{subsection:AssumptionOnSolution}
For any $t$, we assume the exact solution $\psi(\cdot,t)$ to be essentially supported within a cube $\Omega = [-L,L]^N$, 
for given $K$-independent $L$, possibly in a much smaller subregion (as in the case of a moving wavepacket). 
%In such a situation, the Hermite Galerkin approach (possibly with adaptive changes of the with and center of the Hermite basis) 
%appears favorable to a global Galerkin approach on $\Omega$. 
We may define adequately compressed or stretched basis functions
\begin{align*}
 \tilde{\varphi}_{\bfk}(x) = \prod_{l=1}^N\varphi_{k_l}(Sx_l),\qquad \bfk \in \Kc,
\end{align*}
for some positive $S$, to vary resolution. As the univariate $\varphi_K$ is negligibly small outside the interval $[-\sqrt{2(K\!+\!1)}\!-\!1,\sqrt{2(K\!+\!1)}\!+\!1]$, 
we require $S$ and $K$ to be chosen such that $SL \geq \sqrt{2(K+1)}+1$. 
% \begin{align*}
%  SL \geq \sqrt{2(K+1)}+1.
% \end{align*}
For a higher resolution within $\Omega$, we have to increase $S$ and $K$ simultaneously. For ease of presentation, we restrict ourselves to $S = 1$.\\ 
The restrictive assumption on $\psi$ being confined to a fixed cube during propagation in time is indeed dispensable: 
if $\psi(\cdot,t)$ has support outside $\Omega$, a reasonable Galerkin approximation might require a larger choice of $K$, 
or the polynomial approximation of the potential beyond the cube might become useless without increasing $R$. 
This is an issue of  \eqref{align:ODECoefficientsPol} being a good approximation to \eqref{align:SchroedingerEquation}. 
As long as the potential is sufficiently smooth in a region where the solution does not essentially vanish, the fast algorithm is applicable. This is true for the numerical example considered in Section \ref{subsection:timeintegration}. 
A set of moving basis functions 
that adapt to the localisation of $\psi$ might be preferable. 
In fact, the fast algorithm has been applied successfully with a moving wavepacket basis, see \cite{FGL09}. 
However, an adaptation of the following error analysis to this setting complicates the presentation due to the presence of time-dependent evolution parameters in the basis recursion, 
and we restrict our attention to the case of confined $\psi$ and basis functions localized around zero.
%%%%%%%%%%%%%%%%%%%%%%%%%%%%%%%%%%%%%%%%%%%%%%%%%%%%%%%%%%%%%%%%%%%%%%%%%%%%%%%%%%%%%%
%%%%%%%%%%%%%%%%%%%%%%%%%%%%%%%%%%%%%%%%%%%%%%%%%%%%%%%%%%%%%%%%%%%%%%%%%%%%%%%%%%%%%%
\section{Discretization in time}\label{section:SemiDiscretizationInTime}
\subsection{Need for matrix-vector products}\label{subsection:MagnusIntegrators}
The resulting coefficient initial value problem \eqref{align:ODECoefficientsPol} is of the general form $i \dot{y}(t) = A(t) y(t)$ with a time-dependent Hermitian matrix $A(t)$ and initial value $y(0) = y_0$. Polynomial integrators of the form $y^{n+1} = P(-ihA(t))y^n$ as well as discretizations of the matrix 
exponential with a splitting procedure or a Magnus integrator each require multiplications of $A$ with a vector in each time step. In the present paper, we restrict our attention to the latter choice, see  
\cite{BCOR09}, in particular, Sections 5 and 6, for numerical integration methods based on Magnus expansions. Magnus integrators consist of an exponential 
stepping procedure of the form $y^{n+1} = \exp(\Omega^n) y^n$, 
where $y^n \approx y(t_n)$, $t_n = hn$ with time-step size $h$, for a suitable choice of $\Omega^n$. Possible choices are the exponential mid-point rule
\begin{align}\label{align:expMPR}
 \Omega^n = -ihA(t^n + h/2)
\end{align}
or the 2-stage Gauss--Legendre based method with nodes $c_{1,2} = \frac{1}{2} \mp \frac{\sqrt{3}}{6}$,
\begin{align}\label{align:2stageMagnus}
 \Omega^n = -\frac{i}{2}h (A_1 + A_2) - \frac{\sqrt{3}}{12}h^2 [A_2,A_1],\quad A_j = A(t^n + c_j h),\,j = 1,2,
\end{align}
where $[\cdot,\cdot]$ denotes the commutator of matrices. In our setting, we have
\begin{align*}
 A(t) = \mathcal{H}_{\Kc,\text{pol}}(t) = \mathcal{D}_\Kc + \mathcal{W}_{\Kc,\text{pol}}(t).
\end{align*}
\cite{HL03} show that the methods \eqref{align:expMPR} and \eqref{align:2stageMagnus} are of optimal temporal orders 2 and 4, respectively, for the Schr\"odinger equation with a bounded potential. %The schemes \eqref{align:expMPR} and \eqref{align:2stageMagnus} both have a unitary propagator.
%%%%%%%%%%%%%%%%%%%%%%%%%%%%%%%%%%%%%%%%%%%%%%%%%%%%%%%%%%%%%%%%%%%%%%%%%%%%%%%%%%%%%%
\subsection{Lanczos method for the matrix exponential}\label{subsection:LanczosMethod}
We apply the Lanczos method in order to approximate the matrix exponential $\exp(\Omega^n)$, see \cite{L08}, Chapter III.2.2, for a 
more detailed outline including further references and an algorithmic description. Consider a general initial value problem $i \dot{y}(t) = A y(t)$
with an $d \times d$ Hermitian matrix $A$ and $y(0) = y_0$. The Hermitian Lanczos process generates recursively the basis $V_m = (v_1|\dots|v_m) \in \C^{d \times m}$ 
and a tridiagonal coefficient matrix $T_m \in \C^{m \times m}$ such that $T_m = V_m^* A V_m$. 
This requires $m$ multiplications of $A$ on a vector, where $m \ll d$. The matrices $V_m$ and $T_m$ are used to approximate 
\begin{align*}%\label{align:LanczosMethod}
 y(t) = \exp(-itA)y_0 \approx V_m \exp(-itT_m)e_1, \qquad e_1 = (1,0,{\dots},0)^\transpose \in \R^m.
\end{align*}
In our setting, we have $-ihA = \Omega^n$, $y_0 = y^n$. In each time step, 
for all specific choices of $\Omega^n$, this involves the action of $\Wc_{\Kc,\text{pol}}(t)$ on Lanczos vectors $v_k$, 
evaluated at times $t$ depending on the chosen Magnus integrator.
%%%%%%%%%%%%%%%%%%%%%%%%%%%%%%%%%%%%%%%%%%%%%%%%%%%%%%%%%%%%%%%%%%%%%%%%%%%%%%%%%%%%%%
%%%%%%%%%%%%%%%%%%%%%%%%%%%%%%%%%%%%%%%%%%%%%%%%%%%%%%%%%%%%%%%%%%%%%%%%%%%%%%%%%%%%%%
\section{The fast algorithm}\label{section:FastAlgorithm}
%%%%%%%%%%%%%%%%%%%%%%%%%%%%%%%%%%%%%%%%%%%%%%%%%%%%%%%%%%%%%%%%%%%%%%%%%%%%%%%%%%%%%%
\subsection{Coordinate matrices for a direct operation}\label{subsection:FastAlgorithmIdeas}
We consider the product of $\mathcal{W}_{\Kc,\text{pol}}$, see \eqref{align:ODECoefficientsPol}, and an arbitrary vector $v \in \C^{|\Kc|}$. For each direction, we define coordinate matrices over an arbitrary index set $\Kc$ by
\begin{align*}
 X_{\Kc}^{(l)},\quad 1 \leq l \leq N,& & \left(X_{\Kc}^{(l)}\right)_{\bfj\bfk} = (\varphi_\bfj,q_l \varphi_{\bfk}),\quad \bfj,\bfk \in \Kc.
\end{align*}
Due to the orthonormality of the basis and with the help of the one-dimensional recurrence relation \eqref{align:Recurrence1D}, 
the action of  $X_{\Kc}^{(l)}$ on $v$ is given by
\begin{align}\label{align:DirectOperationWithX}\begin{aligned}
 &\left(X_{\Kc}^{(l)}v\right)_\bfj = \sum_{\bfk \in \Kc} (\varphi_\bfj, q_l \varphi_\bfk) v_\bfk\\
  &= \sum_{\bfk \in \Kc} \left(\varphi_\bfj, \sqrt{\frac{k_l+1}{2}}\varphi_{\bfk + \bfe_l} + \sqrt{\frac{k_l}{2}} \varphi_{\bfk-\bfe_l}\right) v_\bfk
  =\!\sqrt{\frac{j_l}{2}}v_{\bfj-\bfe_l}\! +\!\sqrt{\frac{j_l+1}{2}}v_{\bfj+\bfe_l},\end{aligned}
\end{align}
for all  $\bfj \in \Kc$ and $1 \leq l \leq N$. In case $j_l = 0$ or $j_l = k$, the next-to-last or last term vanishes, respectively. 
The matrix-vector product $X_{\Kc}^{(l)}v$ can thus be computed directly in $\Oc(|\Kc|)$ operations.
%%%%%%%%%%%%%%%%%%%%%%%%%%%%%%%%%%%%%%%%%%%%%%%%%%%%%%%%%%%%%%%%%%%%%%%%%%%%%%%%%%%%%%
\subsection{Insertion into the polynomial}\label{subsection:FastAlgorithm}
In case $W^{\text{pol}}(x) = x_l$, for some $l$, the Galerkin matrix reduces to the $l$th coordinate matrix and $\Wc_{\Kc,\text{pol}}v = X_{\Kc}^{(l)}v$. 
As proposed in \cite{FGL09}, the idea is to compute $\Wc_{\Kc,\text{pol}}v$ for any polynomial $W^{\text{pol}}$ by formally inserting the coordinate matrices into the polynomial. 
We can thus compute
\begin{align}\label{align:Horner}\begin{aligned}
 \Wc_{\Kc,\text{pol}}(t)v &\approx W^{\text{pol}}(X_{\Kc},t)v = \sum_{\bfr \in \Rc} \alpha_\bfr(t) \left(\prod_{l=1}^N T_{r_l}\left(\frac{1}{L}X_{\Kc}^{(l)}\right)\right)v \\
  &=  \sum_{\bfr \in \Rc} \alpha_\bfr(t) \left( T_{r_1}\left(\frac{1}{L}X_{\Kc}^{(1)}\right) \cdot \left({\dots}\left(T_{r_N}\left(\frac{1}{L}X_{\Kc}^{(N)}\right)v\right){\dots}\right)\right.,\end{aligned}
\end{align}
%where $T_{r_l}\left(\frac{1}{L}X_{\Kc}^{(l)}\right)v$ is computed recursively with the help of \eqref{align:ChebyshevRecurrence}.
%where the last sum is computed with the help of \eqref{align:ChebyshevRecurrence} using Clenshaw's algorithm, see \cite{C62}, which is done using \eqref{align:DirectOperationWithX}.
where the products of $T_{r_l}\left(\frac{1}{L}X_{\Kc}^{(l)}\right)$ and a vector are computed with the help of the recurrence \eqref{align:ChebyshevRecurrence} using the direct operation \eqref{align:DirectOperationWithX}.
%%%%%%%%%%%%%%%%%%%%%%%%%%%%%%%%%%%%%%%%%%%%%%%%%%%%%%%%%%%%%%%%%%%%%%%%%%%%%%%%%%%%%%
\subsection{Algorithmic description}\label{subsection:AlgorithmicDescription}
The procedures as given in Figure \ref{figure:AlgorithmicDescriptionFastAlgorithm} describe the fast algorithm for the action of 
$\Wc_{\Kc,\text{pol}}(t)$ on a vector $v \in \C^{|\Kc|}$ for given index sets $\Kc(K)$ for the Galerkin basis and $\Rc(R)$ for the polynomial approximation  $W^{\text{pol}}(x,t)$ 
of the potential over a ($K$-independent) cube $\Omega = [-L,L]^N$ with coefficients $\alpha_{\bfr}(t)$.\\ 
%--------------------------------------------------------------------------------------%
\verbdef{\verbfast}{fast_algorithm}
\verbdef{\verbdirect}{direct_op}
\begin{figure}[h!]
\small
\centering
\begin{minipage}[t]{.50\linewidth}
\begin{algorithm}[H]
\SetKwData{input}{input}
\SetKwData{output}{output}
\SetKwData{res}{res}
\SetKwData{temp}{temp}
\hspace{-.2cm}\BlankLine
\hspace{-.2cm}input: $\Kc(K)$, $\Rc(R) \subset \N^N$,\\
\hspace{-.2cm}\hspace{.8cm} $\left(\alpha_{\bfr}(t)\right)_{\bfr\!\in\!\Rc}$, $L\!\in\!\R$, $v\!\in\!\C^{|\Kc|}$\\
\hspace{-.2cm}output: res = $W^{\text{pol}}(X_{\Kc},t)v$\\
\hspace{-.2cm}res $ = 0$\\
\hspace{-.2cm}\For{$\bfr\in\Rc$}{
\hspace{-.2cm}  \For{$l=1$ \KwTo $N$}{
\hspace{-.2cm}    \underline{\textbf{if}$r_l > 0$:}\\
\hspace{-.2cm}    $w_- = v$\\
\hspace{-.2cm}    $w_+ =$\,$\frac{1}{L}$\verbdirect$(\Kc,l,v)$\\
\hspace{-.2cm}    \For{$r=2$ \KwTo $r_l$}{
\hspace{-.2cm}    temp\,$= w_+$\\
\hspace{-.2cm}    $w_+\!=$\,$\frac{2}{L}\,$\verbdirect$(\Kc,l,w_+\!)$\\
\hspace{-.2cm}    $w_+\!= w_+ -w_-$, $w_- =$\, temp\\
\hspace{-.2cm}    }
\hspace{-.2cm}    $v = w_+$\\
\hspace{-.2cm}  }
\hspace{-.2cm}  res $=$ res $+ \alpha_{\bfr}(t)v$
\hspace{-.2cm}}
\caption{\verbfast}
\end{algorithm}
\end{minipage}
\begin{minipage}[t]{.49\linewidth}
\begin{algorithm}[H]
\SetKwData{input}{input}
\SetKwData{output}{output}
\SetKwData{res}{res}
\SetKwData{temp}{temp}
\hspace{-.2cm}\BlankLine
\hspace{-.2cm}input: $\Kc(K) \subset \N^N$,\\
\hspace{-.2cm}\hspace{.8cm} $l \in \{1,{\dots},N\}$, $v \in \C^{|\Kc|}$\\
\hspace{-.2cm}output: res = $X_{\Kc}^{(l)}v$\\
\hspace{-.2cm}\For{$\bfj\in\Kc$}{
\hspace{-.2cm}  res$_{\bfj} = \begin{cases}
	  \sqrt{\frac{1}{2}}v_{\bfj+\bfe_l},& \hspace{-.2cm}j_l\!=\!0,\\
	  \sqrt{\frac{K}{2}}v_{\bfj-\bfe_l},& \hspace{-.2cm}j_l\!=\!K,\\
	  \sqrt{\frac{j_l}{2}}v_{\bfj-\bfe_l} + &\\
	  \hspace{.7cm} \sqrt{\frac{j_l+1}{2}}v_{\bfj+\bfe_l},& \hspace{-.2cm}\text{else}.
	  \end{cases}$
\hspace{-.2cm}}
\caption{\verbdirect}
\end{algorithm}
\end{minipage}
 \caption{Algorithmic description of the fast algorithm \eqref{align:Horner} for a matrix-vector product using \eqref{align:ChebyshevRecurrence} and direct operation \eqref{align:DirectOperationWithX} on a vector with $X_{\Kc}^{(l)}$.}
 \label{figure:AlgorithmicDescriptionFastAlgorithm}
\end{figure}
%--------------------------------------------------------------------------------------%
%--------------------------
Given the time-step size $h$, a number of Lanczos steps $m$ in each time step, 
and the initial coefficient vector $c^0_{\text{pol}}$ of unit norm, time propagation of \eqref{align:ODECoefficientsPol} using a Magnus integrator 
is done as outlined in Figure \ref{figure:AlgorithmicDescriptionTimePropagation}. Step (1) has to be repeated in each time step only in case of a 
time-dependent potential. Step (3) requires a diagonalization of (small) $T_m^{(n)}$, and the product of $V_m^{(n)}$ and a vector costs $\Oc(|\Kc| m^2)$ operations. 
%--------------------------
\verbdef{\verbtime}{time_propagation}
\begin{figure}[h!]
 \centering
 \small
 \begin{minipage}[t]{.98\linewidth}
 \begin{algorithm}[H]
\hspace{-.2cm}\BlankLine
\hspace{-.2cm}input: $\Kc(K), \Rc(R) \subset \N^N$,\, $h,L \in \R$, $m \in \N$,\,$c_{\text{pol}}^0 \in \C^{|\Kc|}$\\
\hspace{-.2cm}output: res = $c_{\text{pol}}^{t_{\text{end}}h^{-1}} \approx c_{\text{pol}}(t_{\text{end}})$\\
\hspace{-.2cm}\For{$n = 0$ \KwTo $t_{\text{end}}h^{-1}$}{
\hspace{-.2cm}  \begin{itemize}
   \item[(1)] Compute coefficients $\alpha_{\bfr}(t)$ of $W^{\text{pol}}(x,t) = \sum_{\bfr\in\Rc} \alpha_{\bfr}(t)T_{\bfr}(x/L)$ as pre-\\scribed by the chosen integrator.
   \item[(2)] Do $m$ Lanczos steps to obtain $V_m^{(n)} = \left(v_1|\dots|v_m\right)$ and $T_m^{(n)}$ starting\\ from $v_1 = c_{\text{pol}}^n$. In each Lanczos step, use:\\
	      \quad $\left(\mathcal{D}_{\Kc}v_k\right)_\bfj = \sum_{l=1}^N (j_l + \frac{1}{2})(v_k)_\bfj$\\
	      \quad $W^{\text{pol}}(X_{\Kc},t)v_k =$ \verbfast$(\Kc,\Rc,(\alpha_{\bfr})_{\bfr\!\in\!\Rc},L,v_k)$
    \item[(3)] Compute $c_{\text{pol}}^{n+1} = V_m^{(n)} \exp(-ih T_m^{(n)}) e_1$.
  \end{itemize}
\hspace{-.2cm}}
\caption{\verbtime}
\end{algorithm}
\end{minipage}
\caption{Time propagation using the fast algorithm with a Magnus integrator.}
\label{figure:AlgorithmicDescriptionTimePropagation}
 \end{figure}
%%%%%%%%%%%%%%%%%%%%%%%%%%%%%%%%%%%%%%%%%%%%%%%%%%%%%%%%%%%%%%%%%%%%%%%%%%%%%%%%%%%%%%
\subsection{Computational complexity}\label{subsection:Complexity}
For fixed $l$, \verb#direct_op# as given in Figure \ref{figure:AlgorithmicDescriptionFastAlgorithm} requires $\Oc(|\Kc|)$ operations. Thus, 
using \eqref{align:ChebyshevRecurrence}, $T_{r_l}\left(\frac{1}{L}X^{(l)}_{\Kc}\right)v$ is computed in $\Oc(|\Kc|\cdot r_l)$ operations, which 
gives $\prod_{l=1}^N T_{r_l}\left(\frac{1}{L}X^{(l)}_{\Kc}\right) v$ in $\Oc(|\bfr|\cdot |\Kc|)$. Therefore, 
\verb#fast_algorithm# allows to compute the product $W^{\text{pol}}(X_{\Kc},t)v$ termwise in
\begin{align*}
 \Oc\left(\sum_{\bfr\in\Rc} |\bfr| \cdot |\Kc|\right)
\end{align*}
operations, i.e., linearly in the size of the basis. Due to $|\Rc| \ll |\Kc|$, the factor $|\Kc|$ is the dominant contribution and, in case of $W$ being time-dependent, the costs for re-computing the interpolation coefficients in each step are negligible. 
In Figure \ref{figure:Complexity}, we compare assembling $\Wc_{\Kc,\text{pol}}$ to a direct computation of $W^{\text{pol}}(X_{\Kc})v$ with respect to CPU time 
for a hyperbolically reduced index set $\Kc$ and a stretched torsional potential %Hénon-Heiles potential 
\begin{align}\label{align:TorsionalPotential}
 W(x) = \sum_{l=1}^N \left(1 - \cos(x_l/L)\right), \quad x \in \Omega, \quad L = 16,%W(x) = \sum_{l=1}^{N-1} (x_l/L)^2(x_{l+1}/L) - \frac{1}{3}(x_{l+1}/L)^3, \quad x \in \Omega,
\end{align}
as approximated by Chebyshev interpolation with $R = 8$ nodes on each coordinate axis (yielding an interpolation error of size $\approx 1\mh10$). The entries of $\Wc_{\Kc,\text{pol}}$ are discretized using Gauss--Hermite quadrature (done as explained in Section \ref{section:GaussHermiteQuadrature}). As the figures reveal, on a hyperbolically reduced index set, the fast algorithm  
lowers the computational effort by several orders of magnitude for reasonable choices of $K$. The larger $K$, the better the reduction (for fixed $N$).  
% --------------------------
 \begin{figure}
 \floatbox[{\capbeside\thisfloatsetup{capbesideposition={right,top},capbesidewidth=6.3cm}}]{figure}[\FBwidth]
 {\caption{Observed CPU times in secs for assembling the hyperbolically indexed matrix $\Wc_{\Kc,\text{pol}}$ and 
 multiplying it with a random vector $v \in \R^{|\Kc|}$ (column (1)) and for the fast algorithm (column (2))  
 with a torsional potential as given in \eqref{align:TorsionalPotential} approximated by Chebyshev 
 interpolation ($R\!=\!8$). Last column: approximate ratios of computation times.
 All figures have been obtained on a desktop computer with an Intel Core 2 
Duo E8400 3.00 GHz processor with 4 GB RAM.
\label{figure:Complexity} }}
 {%\tabulinesep=1.2mm \small
 \begin{tabular} {|l|r|r|r|r|r|}%\begin{tabu}{...}
 \hline $N$& $K$& \multicolumn{1}{c|}{(1)}&\multicolumn{1}{c|}{(2)} &\scriptsize $\approx$(1)/(2)\\
 \hline\hline \multirow{3}{*}{2}	& $20$& \scriptsize $2.83\mh01$			& \scriptsize $1.16\mh03$		& \scriptsize $2.4\mhp02$\\
					& $60$& \scriptsize $\approx\!8.7$ secs	& \scriptsize $1.55\mh03$		& \scriptsize $2.6\mhp03$\\
				       & $100$& \scriptsize $\approx 50$ secs		& \scriptsize $2.87\mh03$		& \scriptsize $1.4\mhp04$\\
 \hline \multirow{3}{*}{3} & $20$& \!\scriptsize $2.51\mhp00$				& \scriptsize $4.57\mh03$		& \scriptsize $5.5\mhp02$\\
			   & $60$& \scriptsize $>\!2.5$ min				& \scriptsize $8.38\mh03$		& \scriptsize $1.5\mhp04$\\
			  & $100$& \scriptsize $> 16$ min				& \scriptsize $1.67\mh02$		& \scriptsize $6.0\mhp04$\\
\hline \multirow{1}{*}{4}& $60$& \scriptsize $\approx 23$ min& \scriptsize $3.29\mh02$& \scriptsize $4.2\mhp04$\\
 \hline
 \end{tabular}}
\end{figure}
% -----------------------------

%%%%%%%%%%%%%%%%%%%%%%%%%%%%%%%%%%%%%%%%%%%%%%%%%%%%%%%%%%%%%%%%%%%%%%%%%%%%%%%%%%%%%%
%%%%%%%%%%%%%%%%%%%%%%%%%%%%%%%%%%%%%%%%%%%%%%%%%%%%%%%%%%%%%%%%%%%%%%%%%%%%%%%%%%%%%%
\section{Relation to Gauss--Hermite quadrature}\label{section:GaussHermiteQuadrature}
\subsection{Product Gauss--Hermite quadrature}\label{subsection:GaussHermiteQuadrature} In preparation for the subsequent error analysis, we consider 
Gaussian quadrature for the weight function $e^{-x^2}$ over $\R$ in each direction (see, e.g., \cite{G12}, Chapter 3.2). 
Let $\xi_m$ denote the zeros of $H_{M+1}$ with corresponding weights $w_m$. 
The resulting quadrature formula $(w_m,\xi_m)_{m=0}^M$ is exact for polynomials  of degree $\leq 2M + 1$. In higher dimensions, we set
\begin{align*}
 \xi_\bfm = (\xi_{m_1},\dots,\xi_{m_N}),& & \omega_\bfm = \prod_{l=1}^N \omega_{m_l} = \prod_{l=1}^N w_{m_l} e^{\xi_{m_l}^2}, & & \bfm \in \mathcal{M}_{\text{full}}(M),
\end{align*}
with a full $N$-dimensional index set $\mathcal{M}_{\text{full}}(M)$. This yields a product quadrature 
\begin{align}\label{eq:GHQuad_full}
\begin{aligned}
 (\mathcal{W}_{\Kc,\text{pol}}(t))_{\bfj\bfk} &\approx (\mathcal{W}_{\Kc,\text{pol}}^{\text{GH}(M)}(t))_{\bfj\bfk} = \sum_{\bfm \in \mathcal{M}} \omega_\bfm \varphi_\bfj(\xi_\bfm)W^{\text{pol}}(\xi_\bfm,t)\varphi_\bfk(\xi_\bfm)\\
  &= \sum_{\bfr\in\Rc} \alpha_{\bfr}(t) \prod_{l=1}^N \sum_{m_l=0}^M \underbrace{\omega_{m_l}\varphi_{j_l}(\xi_{m_l})T_{r_l}(\xi_{m_l}/L)\varphi_{k_l}(\xi_{m_l})}_{\star},
\end{aligned}
\end{align}
which is exact if $W^{\text{pol}}(\cdot,t)H_{\bfj}H_{\bfk}$ is a polynomial of degree $\leq 2M + 1$ in each direction, 
where $H_\bfk(x) = \prod_{l=1}^N H_{k_l}(x_l)$ is a tensor product of univariate Hermite polynomials. 
Having obtained $\xi_m$ and $w_m$ (see \cite{NR07}, Chapter 4.6), the recursions \eqref{align:Recurrence1D} and \eqref{align:ChebyshevRecurrence} 
allow us to compute $\varphi_k(\xi_m)$ and $T_r(\xi_m)$, for $0 \leq k \leq K$, $0 \leq r \leq R$, $0 \leq m \leq M$, in $\Oc((R+K)M)$ in advance. 
Thus, given the terms $\star$, we compute the whole expression \eqref{eq:GHQuad_full} in $\Oc(|\Rc|MN)$, 
and the assembly of $\mathcal{W}_{\Kc,\text{pol}}^{\text{GH}(M)}$ requires $\Oc\left(|\Kc|^2 \cdot |\Rc| \cdot M \cdot N\right)$ operations. We choose $M=K$, see Lemma  \ref{Lemma:QuadratureEqualsInsertion}.\\
Smolyak sparse grid quadrature (see \cite{S63, Z91, GG98}) adapted to the 
increasingly oscillatory behavior of the high-order Hermite functions is discussed in \cite{L08}, Chapter III.1.2., 
where it is pointed out that a sufficiently accurate sparse grid 
quadrature requires at least $\Oc(|\Kc|^2 \cdot M)$ evaluations of the potential anyway.
%%%%%%%%%%%%%%%%%%%%%%%%%%%%%%%%%%%%%%%%%%%%%%%%%%%%%%%%%%%%%%%%%%%%%%%%%%%%%%%%%%%%%%
\subsection{Insertion of coordinate matrices revisited}
We define $X^{(l)}_{\text{full}} = X^{(l)}_{\Kc_{\text{full}}}$ to be the coordinate matrices over the full index set. The matrices 
\begin{align*}
 \Xi^{(l)} = \text{diag}_{\bfk \in \Kc} (\xi_{k_l}) \in \R^{|\Kc_{\text{full}}|\times|\Kc_{\text{full}}|},& & U_{\bfj\bfk} = \sqrt{\omega_\bfj}\varphi_\bfk(\xi_\bfj),\quad \bfj, \bfk \in \Kc_{\text{full}},
\end{align*}
yield a diagonalization $X^{(l)}_{\text{full}}=U^\transpose \Xi^{(l)} U\in\R^{|\Kc_{\text{full}}|\times|\Kc_{\text{full}}|}$, which is readily seen from 
\begin{align*}%\label{align:CheckDiagonalizationND}
\begin{aligned}
 (U^\transpose \Xi^{(l)} U)_{\bfj\bfk} &= \sum_{\bfm \in \Kc} \omega_\bfm \xi_{m_l} \varphi_\bfj(\xi_\bfm) \varphi_\bfk (\xi_\bfm) = \left(\varphi_\bfj, q_l \varphi_\bfk\right)^{\text{GH}(K)} = (X^{(l)}_{\text{full}})_{\bfj\bfk},
\end{aligned}
\end{align*}
by the fact that there are exactly $K+1$ quadrature nodes in each direction and that this yields an exact integration. The matrix $U$ is unitary, which follows from orthonormality of the basis and
\begin{align*}
 (U^\transpose U)_{\bfj\bfk} = \sum_{\bfm \in \Kc} u_{\bfm\bfj}u_{\bfm\bfk} = \sum_{\bfm \in \Kc} \omega_\bfm \varphi_\bfj(\xi_\bfm)\varphi_\bfk(\xi_\bfm) = (\varphi_\bfj, \varphi_\bfk)^{\text{GH}(K)} = (\varphi_\bfj, \varphi_\bfk) = \delta_{\bfj\bfk}.
\end{align*}
This allows to compute
\begin{align}\label{align:Xhochr}
 X^\bfr_{\text{full}}=\left(X^{(1)}_{\text{full}}\right)^{r_1}\!\!\!{\dots}\left(X^{(N)}_{\text{full}}\right)^{r_N}\!\!\!
        =\prod_{l=1}^N\left(U^\transpose \text{diag}\!\left(\xi_{m_l}^{r_l}\right)U\right)=U^\transpose \text{diag}\!\left(\xi_\bfm^\bfr\right) U,
\end{align} 
and we get the following
\begin{lemma}\label{Lemma:QuadratureEqualsInsertion}
Choosing $\mathcal{M} = \Kc_{\textnormal{full}}$ (i.e., $M = K$) for the full product quadrature and basis index sets, respectively, we get
\begin{align*}
 W^{\text{pol}}(X_{\textnormal{full}},t)_{\bfj\bfk} = (\mathcal{W}_{\Kc_{\textnormal{full}},\textnormal{pol}}^{\textnormal{GH}(K)}(t))_{\bfj\bfk},\quad \bfj,\bfk \in \Kc_{\textnormal{full}},
\end{align*}
where $W^{\textnormal{pol}}(X_{\textnormal{full}},t)$ denotes formal insertion of $X_{\textnormal{full}}^{(l)}$ into $W^{\textnormal{pol}}$ according to \eqref{align:Xhochr}.\hfill\dollar\\
\end{lemma}
This argument is common in the context of DVR techniques, see \cite{LC00}. 
The ordering of the factors $\left(X^{(l)}_{\text{full}}\right)^{r_l}$ in $W^{\text{pol}}(X_{\text{full}},t)$ is arbitrary.\\
Deriving the equivalence of full product quadrature and formal insertion requires a bijection $\mathcal{M} \leftrightarrow \Kc_{\text{full}}$. 
Simultaneously reducing $\mathcal{M}$ 
and $\Kc_{\text{full}}$ invalidates the exactness of the Gauss--Hermite quadrature, reducing only $\Kc_{\text{full}}$ makes the above diagonalization argument
no longer correct at all. For a 
reduced index set $\Kc \subsetneq \Kc_{\text{full}}$, an assertion analogous to Lemma 
\ref{Lemma:QuadratureEqualsInsertion} can therefore not be expected. In the fast algorithm, 
we employ the above reduced coordinate matrices $X^{(l)}_{\Kc}$. Hence, we expect the fast algorithm to induce errors due to quadrature and index set reduction, and a consideration of full product quadrature facilitates 
the error analysis.\\
%%%%%%%%%%%%%%%%%%%%%%%%%%%%%%%%%%%%%%%%%%%%%%%%%%%%%%%%%%%%%%%%%%%%%%%%%%%%%%%%%%%%%%
%%%%%%%%%%%%%%%%%%%%%%%%%%%%%%%%%%%%%%%%%%%%%%%%%%%%%%%%%%%%%%%%%%%%%%%%%%%%%%%%%%%%%%
\section{Error analysis}\label{section:ErrorAnalysis}
%%%%%%%%%%%%%%%%%%%%%%%%%%%%%%%%%%%%%%%%%%%%%%%%%%%%%%%%%%%%%%%%%%%%%%%%%%%%%%%%%%%%%%
\subsection{Preliminaries}\label{subsection:ErrorPreliminaries}
\textit{Definition of errors:} Consider an arbitrary vector $v \in \C^{|\Kc|}$. We are interested in computing the product $\Wc_{\Kc,\text{pol}}(t)v$ with a matrix $\Wc_{\Kc,\text{pol}}$ 
as given in Section \ref{subsection:ApproximationPotential}. The fast algorithm as developed in Section \ref{section:FastAlgorithm} gives rise to an error due to quadrature and to an error due to index set reduction, the former being given by
\begin{align}\label{align:ErrorQuadratureDef}
 E^{\text{quad}} = \left(E_{\bfj,\bfk}\right)_{\bfj,\bfk \in \Kc},& & E_{\bfj,\bfk} =  (\Wc_{\Kc,\text{pol}}(t))_{\bfj\bfk} - (\Wc_{\Kc,\text{pol}}^{\text{GH}(K)}(t))_{\bfj\bfk}.
\end{align}
Formally inserting the reduced coordinate matrices into the polynomial yields an error
\begin{align*}
W^{\text{pol}}(X_{\Kc},t)v - \Wc_{\Kc,\text{pol}}^{\text{GH}(K)}(t)v&= \left[W^{\text{pol}}(X_{\Kc},t)v - \Omega(W^{\text{pol}}(X_{\text{full}},t))v\right]\\
  &\hspace{2.3cm}+ \underbrace{\left[\Omega(W^{\text{pol}}(X_{\text{full}},t)) - \Wc_{\Kc,\text{pol}}^{\text{GH}(K)}(t)\right]}_{\star}v,
\end{align*}
where the operator
\begin{align}\label{align:DefOmegas}
 \Omega: \C^{|\Kc_{\text{full}}|\times|\Kc_{\text{full}}|} \to \C^{|\Kc| \times |\Kc|},\quad \Omega(A) = (A_{\bfj\bfk})_{\bfj,\bfk \in \Kc}
\end{align}
cuts a fully indexed matrix to a reduced index set. The difference $\star$ vanishes by virtue of Lemma \ref{Lemma:QuadratureEqualsInsertion}. One easily verifies 
\begin{align*}
 \Omega(W^{\text{pol}}(X_{\text{full}},t))v = \Omega\left(W^{\text{pol}}(X_{\text{full}},t)\Omega_+(v)\right)
\end{align*}
where the operator
\begin{align*}
\Omega_+\!: \C^{|\Kc|}\!\!\to\!\C^{|\Kc_{\text{full}}|},\qquad\left(\Omega_+(v)\right)_{\bfj}=\begin{cases} v_\bfj,&\bfj \in \Kc,\\ 0,&\bfj \notin \Kc,\end{cases} 
\end{align*}
blows up a vector with zeros at indices missing in $\Kc$, and $\Omega$ is defined as in \eqref{align:DefOmegas} both for matrices and vectors. Hence, the error due to index set reduction is given by
\begin{align}\label{align:ErrorGridReductionDef}
E^{\text{red}}(v) = \left(E_\bfj\right)_{\bfj \in \Kc},\,\, E_\bfj = \left(W^{\text{pol}}(X_{\Kc},t)v\right)_{\bfj} - \left(W^{\text{pol}}(X_{\text{full}},t)\Omega_+(v)\right)_\bfj,\quad \bfj \in \Kc.
\end{align}
\textit{Assumption:} For the following error analysis, we make the general decay assumption
\begin{align}\label{align:GeneralDecayCondition}
 v_{\bfk} = \Oc\left(\prod_{l=1}^N \max(k_l,1)^{-\beta}\right),\quad\quad \bfk \in \Kc,
\end{align}
for the vector coefficients of $v$, with some $\beta \in \N$. Thus, the larger its index, the smaller the vector component. Assumption \eqref{align:GeneralDecayCondition} reflects the natural decay behavior of the coefficients 
in a product Hermite expansion of a sufficiently smooth function over a hyperbolically reduced index set, see \cite{L08}, Thm. III.1.5. 
It is used in Sections \ref{subsection:ErrorQuadrature} and \ref{subsection:ErrorGridReduction} to 
compensate large error components in matrix-vector products.\\
Finally, for $\bfr \in \N^N$, we set $r_{\max}(\bfr) = \max_{1 \leq l \leq N} r_l$.
%%%%%%%%%%%%%%%%%%%%%%%%%%%%%%%%%%%%%%%%%%%%%%%%%%%%%%%%%%%%%%%%%%%%%%%%%%%%%%%%%%%%%%
\subsection{Error $E^{\text{quad}}$ due to quadrature}\label{subsection:ErrorQuadrature}
%----------------------------------
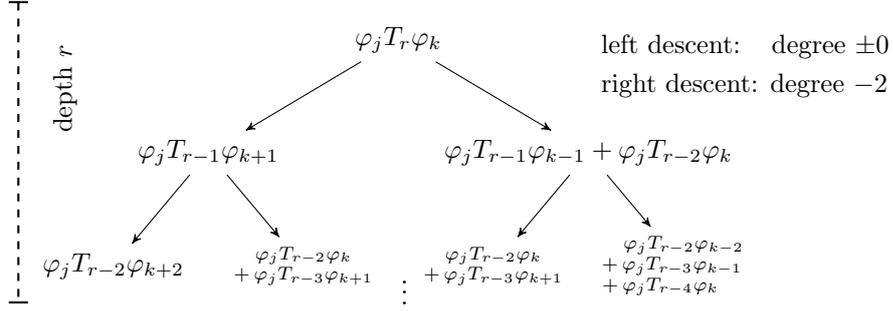
\begin{figure}[t!]
 \begin{tikzpicture}[->,>=stealth',level/.style={sibling distance = 5cm/#1,level distance = 1.5cm}] 
\node [] {$\varphi_jT_r \varphi_k$}
    child{ node [] {$\varphi_jT_{r-1}\varphi_{k+1}$} 
            child{ node [] {$\varphi_jT_{r-2}\varphi_{k+2}$} edge from parent node[left] {}}
            child{ node [] {$\substack{\varphi_jT_{r-2}\varphi_{k}\\+\,\varphi_jT_{r-3}\varphi_{k+1}}$} edge from parent node[right] {}}
	edge from parent node[left,pos=.2] {}
    }
    child{ node [] {$\varphi_jT_{r-1}\varphi_{k-1} + \varphi_jT_{r-2}\varphi_k$}
            child{ node [] {$\substack{\varphi_jT_{r-2}\varphi_{k}\\+\,\varphi_jT_{r-3}\varphi_{k+1}}$} edge from parent node[left] {}}
            child{ node [] {$\substack{\varphi_jT_{r-2}\varphi_{k-2}\\\hspace{-.3cm}+\,\varphi_jT_{r-3}\varphi_{k-1}\\\hspace{-.3cm}+\,\varphi_jT_{r-4}\varphi_{k\,\,\,\,\,\,\,}}$} edge from parent node[right] {}}
	edge from parent node[right,pos=.2] {}
    }
;
 \draw [thick,dashed,|-|] (-5,.51) -- (-5,-3.5);
\end{tikzpicture}
\begin{picture}(5,5)
 \put(-60,95){left descent: \hspace{.1cm} degree $\pm 0$}
 \put(-60,80){right descent: degree $-2$}
 \put(-135,0){$\vdots$}
 \put(-265,65){\rotatebox{90}{depth $r$}}
\end{picture}
\caption{Expansion of 1D quadrature error as a binary tree. In case $r\!=\!1$, we just use \eqref{align:Recurrence1D} on a term, in case $r\!=\!0$, 
terms are added to the left child without expanding.}
\label{figure:ErrQuadBinaryTree}
\vspace{-.2cm}
\end{figure}
%----------------------------------
%\begin{Thm}\label{Theorem:ErrorQuadratureMatrixVector}
\begin{theorem} \label{Theorem:ErrorQuadratureMatrixVector}
Let $W^{\textnormal{pol}}(\cdot,t) \approx W(\cdot,t)$ 
be the Chebyshev interpolation polynomial of the potential $W$ on $\Omega = [-L,L]^N$ over $\Rc(R)$ for fixed $L$. Let $\Kc(K)$ be 
a hyperbolically reduced index set with $K \gg R$. 
Then, under assumption \eqref{align:GeneralDecayCondition} on $v \in \C^{|\Kc|}$ (i.e., componentwise decay of order $\beta \in \N$), 
the error due to quadrature behaves as
\begin{align*}
 \left|\left(\left(\Wc_{\Kc,\textnormal{pol}}(t)-\Wc_{\Kc,\textnormal{pol}}^{\textnormal{GH}(K)}(t)\right)v\right)_\bfj\right| \leq C(\Rc,W,L,t) K^{-\beta},\quad\quad\bfj \in \Kc,
\end{align*}
where the matrices $\Wc_{\Kc,\textnormal{pol}}(t)$ and $\Wc_{\Kc,\textnormal{pol}}^{\textnormal{GH}(K)}(t)$ are defined according to Sections \ref{subsection:ApproximationPotential} and 
\ref{subsection:GaussHermiteQuadrature}, respectively. The constant $C(\Rc,W,L,t)$ behaves as in 
\eqref{align:CQuad}, see below, and depends only on $\Rc$, the regularity of $W$, $L$, and time $t$.%\hfill\dollar
\end{theorem}
%\end{Thm}
\emph{Proof.} Termwise consideration of $W^{\text{pol}}$ gives rise to error matrices
\begin{align*}%\label{align:ErrorQuadratureDecomposition}
  E^{\text{quad}}_{\bfj\bfk}\!\!=\!\!\sum_{\bfr \in \Rc}\!\alpha_{\bfr}(t) E^{\bfr}_{\bfj\bfk},\,\, E^{\bfr}_{\bfj\bfk}\!=\!\left(\varphi_{\bfj}(S\cdot),T_{\bfr}(\cdot/L)\varphi_{\bfk}(S\cdot)\right)\!-\!\left(\varphi_{\bfj}(S\cdot),T_{\bfr}(\cdot/L)\varphi_{\bfk}(S\cdot)\right)^{\text{GH}(K)}. 
\end{align*}
\textit{Conversion of 1D-error into binary tree:} In one dimension, applying recursions \eqref{align:ChebyshevRecurrence} and \eqref{align:Recurrence1D} yields a decomposition
\begin{align}\label{align:QuadratureRecurrence1D}\begin{aligned}
 \varphi_j(Sx)T_r(x/L)\varphi_k(Sx) = & \left(\frac{\sqrt{2(k+1)}}{SL}\varphi_j(Sx)T_{r-1}(x/L)\varphi_{k+1}(Sx)\right)\\ 
				       & \hspace{-3cm}+ \left(\frac{\sqrt{2k}}{SL}\varphi_j(Sx)T_{r-1}(x/L)\varphi_{k-1}(Sx) - \varphi_j(Sx)T_{r-2}(x/L)\varphi_k(Sx)\right).
				       \end{aligned}
\end{align}
Due to $SL \geq \sqrt{2(K+1)}+1$ (see Section \ref{subsection:ApproximationPotential}), all coefficients are bounded by
\begin{align}\label{align:QuadratureCoefficientsBound}
 \frac{\sqrt{2k}}{SL} \leq \frac{\sqrt{2(k+1)}}{SL} \leq 1.
\end{align}
Termwise $r$-fold application of \eqref{align:QuadratureRecurrence1D} yields the binary tree pattern $\T$ as given in Figure \ref{figure:ErrQuadBinaryTree} 
(arguments and coefficients omitted). %%%, see Figure \ref{figure:ErrQuadBinaryTree2}. 
Descending right reduces the polynomial degree by $2$, descending left leaves it unaltered. In case $r > k$, 
we may define $\varphi_k(x) = 0$ for $k \leq -1$, preserving the recurrence relation \eqref{align:Recurrence1D} for negative indices. We expand termwise until each leaf carries a single term of the form 
$\varphi_j \varphi_{k+\lambda-\rho}$, where $\rho$ and $\lambda$ are the numbers of index-changing right and left descents, respectively, and $\lambda + \rho \leq r$. By the same procedure for the corresponding quadrature formulas, we convert
$E_{jk}^r$ into a binary tree $\T$ of depth $r$. We examine non-vanishing leaves in $\T$, how many these are, and what quantity they sum up to.\\
\textit{Characterization of non-vanishing leaves:} Fix a leaf in $\T$ connected to the root $E^r_{jk}$ by $\rho$ and $\lambda$ right and left descents, respectively.  For the quadrature error not to vanish, 
we require $2(K+1) \leq j + k + \lambda - \rho$ for each single term. For the exact integral not to vanish, orthogonality yields the requirement  $j = k + \lambda - \rho$, leading to the contradiction $2(K+1) \leq 2j$. Thus, leaves with non-vanishing quadrature errors carry only the quadrature formulas.\\
\textit{Number of non-vanishing leaves:} The condition for a non-vanishing quadrature error at a particular leaf yields $2(K+1) \leq j + k + \lambda -\rho \leq j + k + r - 2\rho$. We define
\begin{align*}
 \rho_{\max}(j,k,r) = \left\lfloor \frac{j+k+r-2(K+1)}{2} \right\rfloor_+ \leq \frac{r}{2}-1,
\end{align*}
which is the maximal number of right descents that does not reduce the polynomial degree of the integrand sufficiently for exact quadrature, where 
$a_+ = \max(a,0)$. In an arbitrary full binary tree of depth $r$, the number of leaves connected to the root by a path containing 
exactly $s$ right descents equals $\binom{r}{s}$. Hence, the number of non-vanishing leaves in $\T$ is given by $a(j,k,r) = \sum_{s=0}^{\rho_{\max}} \binom{r}{s}$. To investigate further, consider a sum $\sum_{s=a}^b \binom{c}{s}$ with $a, b, c \in \N$ and $a \leq 2b < c$. For $s \leq b$,
\begin{align*}
 \binom{c}{s} \left/ \binom{c}{b}\right. = \frac{b! (c-b)!}{s! (c-s)!} = \frac{b\cdot{\dots}\cdot (s+1)}{(c-s)\cdot{\dots}\cdot (c-b+1)} < \left(\frac{b}{c-b+1}\right)^{b-s}.
\end{align*}
The assumption $c > 2b$ implies $b/(c-b+1) < 1$. Therefore,
\begin{align}\label{align:SumBinomialCoeff}
\sum_{s=a}^b\!\binom{c}{s}\!=\!\binom{c}{b}\!\sum_{s=a}^b\!\binom{c}{s}\!\left/\!\binom{c}{b}\right.\!<\!\binom{c}{b}\!\sum_{s=a}^b\!\left(\frac{b}{c-b+1}\right)^{b-s}\!\!\!<\!\binom{c}{b} (b-a+1).
%  &< \binom{c}{b} \sum_{s=0}^{b-a} \left(\frac{2b}{c}\right)^{s} < \binom{c}{b} \left(1-\frac{2b}{c}\right)^{-1} \left(1-\left(\frac{2b}{c}\right)^{b-a+1}\right).\end{aligned}
% \begin{aligned}
%  &\sum_{s=a}^b\!\binom{c}{s}\!=\!\binom{c}{b}\!\sum_{s=a}^b\!\binom{c}{s}\!\left/\!\binom{c}{b}\right.\!<\!\binom{c}{b}\!\sum_{s=a}^b\!\left(\frac{b}{c-b+1}\right)^{b-s}\!\!\!\!\!=\!\binom{c}{b}\!\sum_{s=0}^{b-a}\!\left(\frac{b}{c-b+1}\right)^{s}\\
%  &< \binom{c}{b} \sum_{s=0}^{b-a} \left(\frac{2b}{c}\right)^{s} < \binom{c}{b} \left(1-\frac{2b}{c}\right)^{-1} \left(1-\left(\frac{2b}{c}\right)^{b-a+1}\right).\end{aligned}
\end{align}
By definition, we have $2 \rho_{\max} < r$, thus, 
\begin{align}\label{align:QuadratureError1D_a}
 a(j,k,r)\!<\!\binom{r}{ \left\lfloor\frac{r}{2}-1\right\rfloor_+}\!\frac{r}{2}.
\end{align}
 At best, $j + k + r = 2K + 2$, and we have $a(j,k,r) = 1$. At worst, $j, k = K$, thus, $\rho_{\max} = \left\lfloor\frac{r}{2}-1\right\rfloor_+$, which makes the last estimate almost sharp.\\
\textit{Error accumulation in 1D:} Taking into account boundedness of the coefficients \eqref{align:QuadratureCoefficientsBound} and vanishing exact integrals at non-vanishing leaves, summing up yields
\begin{align}\label{align:QuadratureError1D_b}
\begin{aligned}
 \left|E_{jk}^r\right| \leq a(j,k,r) \cdot \hspace{-.4cm}\max_{\substack{0 \leq j \leq K, -r \leq k \leq K+r\\ j+k \geq 2K + 2}} \left|\sum_{m=0}^K\! \omega_m \varphi_j(\xi_m) \varphi_{k}(\xi_m)\right| = a(j,k,r) \cdot \mu(K,r).\end{aligned}
\end{align}
Due to cancellation effects by Hermite function evaluations with rapidly 
alternating signs, the term $\mu(K,r)$ is of size $\Oc(1)$.\\
\textit{Decomposition of error in multiple dimensions:} We set $\mathcal{N} = \{1,\dots,N\}$ and consider the error matrix for $N \geq 2$. 
For arbitrary $\bfj, \bfk \in \Kc$ and $\bfr \in \Rc$: If $k_l \leq K - r_l + 1$, for all $l \in \Nc$, we have $j_l + k_l + r_l \leq 2K + 1$, and the one-dimensional 
error matrices $E^{r_l}_{j_lk_l}$ vanish. Hence, $E^{\bfr}_{\bfj\bfk} = 0$. Conversely, for fixed $\bfj \in \Kc$ and $\bfr \in \Rc$, 
if there is $\bfk \in \Kc$ such that $E^{\bfr}_{\bfj\bfk} \neq 0$, we find a subset of components $\Nct = \Nct(\bfk) \subseteq \Nc$ such that, for every $l \in \Nct$, $k_l \geq K - r_l + 2$, and $E_{j_l k_l}^{r_l}$ 
does not vanish. This allows for a decomposition (omitting factors $M$ and $L^{-1}$ in $\varphi$ and $T$, respectively)
\begin{align}\label{align:QuadratureErrorDecomposition}\begin{aligned}
 E_{\bfj\bfk}^{\bfr} &= (\varphi_\bfj,T_{\bfr}\varphi_\bfk) - (\varphi_\bfj,T_{\bfr}\varphi_\bfk)^{\text{GH}(K)} = \left[\prod_{l=1}^N (\varphi_{j_l},T_{r_l}\varphi_{k_l})\right] - \left[\prod_{l=1}^N (\varphi_{j_l},T_{r_l}\varphi_{k_l})^{\text{GH}(K)}\right]\\
    &= \left\{\underbrace{\left[\prod_{l \in \tilde{N}} (\varphi_{j_l},T_{r_l}\varphi_{k_l})\right]}_{\textbf{A}}-\underbrace{\left[\prod_{l \in \tilde{N}} (\varphi_{j_l},T_{r_l}\varphi_{k_l})^{\text{GH}(K)}\right]}_{\textbf{B}}\right\}\underbrace{\left[\prod_{l \notin \tilde{N}} (\varphi_{j_l},T_{r_l}\varphi_{k_l})\right]}_{\textbf{C}}
    \end{aligned}
\end{align}
of a non-vanishing entry $E_{\bfj\bfk}^{\bfr}$. On a hyperbolically reduced set $\Kc$, non-vanishing errors $E_{\bfj\bfk}^{\bfr}$ have indices $\bfk$ satisfying
\begin{align*}%\label{align:ConditionHypCross}
 K\!+\!1\!\geq\!\prod_{l\in\Nc} (k_l\!+\!1)\!=\!\left(\prod_{l \notin \Nct} (k_l\!+\!1)\!\right)\!\!\left(\prod_{l \in \Nct} (k_l\!+\!1)\!\right)\!\geq\!\left(\prod_{l \notin \Nct} (k_l\!+\!1)\!\right)\!\!\left(\prod_{l \in \Nct} (K - r_l\!+\!3)\!\right).
\end{align*}
Clearly, for every $\bfk \in \Kc$, if $K \gg r_{\max}$, then $|\tilde{\mathcal{N}}(\bfk)| \leq 1$. Thus, the terms $\textbf{A}$ and $\textbf{B}$ consist 
of exactly one factor each, and $\textbf{A} - \textbf{B}$ equals the one-dimensional quadrature error $E_{j_{l_0}k_{l_0}}^{r_{l_0}}$ for some $l_0 \in \Nc$.\\
\textit{Error estimation in multiple dimensions:} Consider a non-vanishing entry $E^{\bfr}_{\bfj\bfk}$. By the above considerations for the one-dimensional case, 
the term $\textbf{A}$ vanishes. Due to $|\Nct(\bfk)| = 1$, there is $l_0 \in \Nc$ such that $\textbf{B}$ equals $E_{j_{l_0}k_{l_0}}^{r_{l_0}}$. 
For a factor in $\textbf{C}$, using Cauchy--Schwarz, we find $\left(\varphi_{j_l},T_{r_l}\varphi_{k_l}\right) \leq C(r_l)$. Thus, from \eqref{align:QuadratureError1D_a} and \eqref{align:QuadratureError1D_b}, we have
\begin{align*}
 E_{\bfj\bfk}^{\bfr} = \Oc\left(E_{j_{l_0}k_{l_0}}^{r_{l_0}}\right) = \Oc\left(\binom{r_{l_0}}{\left\lfloor \frac{r_{l_0}}{2}-1 \right\rfloor_+}\frac{r_{l_0}}{2}\right).
\end{align*}
For a hyperbolically reduced index set $\Kc(K)$ with $R \leq \frac{K-1}{2}+2$, it is easily seen that the total number of non-vanishing entries in $E^{\bfr}$ is at most $\frac{1}{2}r_{\max}(r_{\max}-1)$. Multiplying the matrix with 
a rapidly decaying vector, we thus find
\begin{align}\label{align:QuadratureErrorNDEstimate}\begin{aligned}
 \left(E^{\bfr} v\right)_\bfj = \Oc\left(\frac{r_{\max}^2(r_{\max}-1)}{4}\binom{r_{\max}}{\left\lfloor \frac{r_{\max}}{2}-1 \right\rfloor_+}(K-r_{\max}+2)^{-\beta}\right).\end{aligned}
\end{align}
Summing up, the time-dependent and rapidly decaying, $L$-dependent interpolation coefficients enter into the constant
\begin{align}\label{align:CQuad}
 C(\Rc,W,L,t) = \sum_{\bfr \in \Rc} \alpha_{\bfr}(t)\frac{1}{4}r_{\max}^2(r_{\max}-1) \binom{r_{\max}}{\left\lfloor \frac{r_{\max}}{2}-1 \right\rfloor_+},%= \Oc\left(\sum_{\bfr \in \Rc}r_{\max}^{3-p}\binom{r_{\max}}{\lfloor r_{\max}/2 \rfloor}2^{-r_{\max}}\right),
\end{align}
and, together with the assumption $K \gg R$, this proves the claim.\hfill\dollar
%%%%%%%%%%%%%%%%%%%%%%%%%%%%%%%%%%%%%%%%%%%%%%%%%%%%%%%%%%%%%%
\subsection{Error $E^{\text{red}}$ due to index set reduction}\label{subsection:ErrorGridReduction}
% ------------------------
\begin{figure}
\tikzstyle{level 1}=[level distance=2.3cm, sibling distance=2.7cm]
\tikzstyle{level 2}=[level distance=3.2cm, sibling distance=1.6cm]
\begin{tikzpicture}[grow=right, sloped,thick,->]
\node[] {$T_{r_l,\bfj}$}
    child {
        % UNTEN
        node[] {\hspace{.8cm} $\sum\limits_{\substack{s=1,\\s\text{ odd}}}^{r_l-1}(-1)^{s+1}T_{r_l-s, \bfj-\bfe_l} \pm v_{,\bfj}$}        
	    child {
	        % UNTEN-UNTEN
                node[label=right:
                    {\hspace{-.8cm} $\substack{\sum\limits_{\substack{s=1,\\s\text{ odd}}}^{r_l-1}\sum\limits_{\substack{t=1,\\t\text{ odd}}}^{r_l-s-1}(-1)^{s+t}T_{r_l-s-t, \bfj-2\bfe_l}\\ \hspace{2.3cm} + \sum\limits_{\substack{s=1,\\s\text{ odd}}}^{r_l-1} (-1)^{s+1}v_{\bfj-2\bfe_l}\,\pm\,v_{\bfj}}$}] {}
                edge from parent
                node[above] {}
                node[below] {}%{\quad$\frac{\sqrt{2(j_l-1)}}{L}$}
            }
            child {
                % UNTEN-OBEN
                node[label=right:
                    {$\left(\dots\right)_{\bfj} + \sum\limits_{\substack{s=1,\\s\text{ odd}}}^{r_l-1} (-1)^{s+1}v_{\bfj}$}] {}
                edge from parent
                node[above] {}%{$\frac{\sqrt{2j_l}}{L}$}
                node[below] {}
            }
        edge from parent         
            node[above] {}
            node[below] {}%{$\frac{\sqrt{2j_l}}{L}$}
    }
    % OBEN
    child {
        node[] {$\sum\limits_{\substack{s=1,\\s\text{ odd}}}^{r_l-1}(-1)^{s+1}T_{r_l-s, \bfj+\bfe_l}$}        
            child {
		% OBEN-UNTEN
                node[ label=right:
                    {$\left(\dots\right)_{\bfj} + \sum\limits_{\substack{s=1,\\s\text{ odd}}}^{r_l-1} (-1)^{s+1}v_{\bfj}$}] {}
                edge from parent
                node[above] {}
                node[below] {}%{\quad$\frac{\sqrt{2(j_l+1)}}{L}$}
            }
            child {
		% OBEN-OBEN
                node[ label=right:
                    {$\left(\dots\right)_{\bfj+2\bfe_l} + \sum\limits_{\substack{s=1,\\s\text{ odd}}}^{r_l-1} (-1)^{s+1}v_{\bfj+2\bfe_l}$}] {}
                edge from parent
                node[above] {}%{$\frac{\sqrt{2(j_l+2)}}{L}$}
                node[below] {}
            }
            edge from parent 
            node[above] {}%{$\frac{\sqrt{2(j_l+1)}}{L}$}
            node[below] {}
     };
 \end{tikzpicture}
\begin{picture}(5,5)
 \put(-45,140){left descent: \hspace{.1cm} $-1$}
 \put(-45,125){right descent: $+1$}
 \put(-35,95){$\dots$}
 %\put(-210,125){\rotatebox{90}{depth $r_l$}}
\end{picture}
\caption{Expansion of $T_{r_l,\bfj}$ as a binary tree $\T^{(l)}(\bfj)$ of depth $r_l$ (90 ° rotation, coefficients omitted). 
The figure illustrates the case of $r_l$ being even. Terms with reduced or unaltered index are attached to the left child, otherwise to the right child.}
\label{figure:ErrRedBinaryTree}
\vspace{-.2cm}
\end{figure}
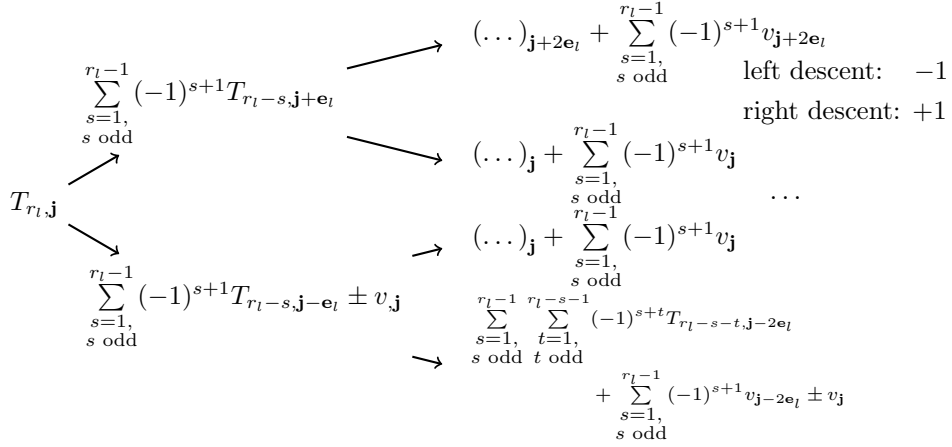
\begin{wrapfigure}{R}{0.5\textwidth}
\includegraphics[trim=5.5cm 14.2cm 5.2cm 11.2cm,clip,scale=.63]{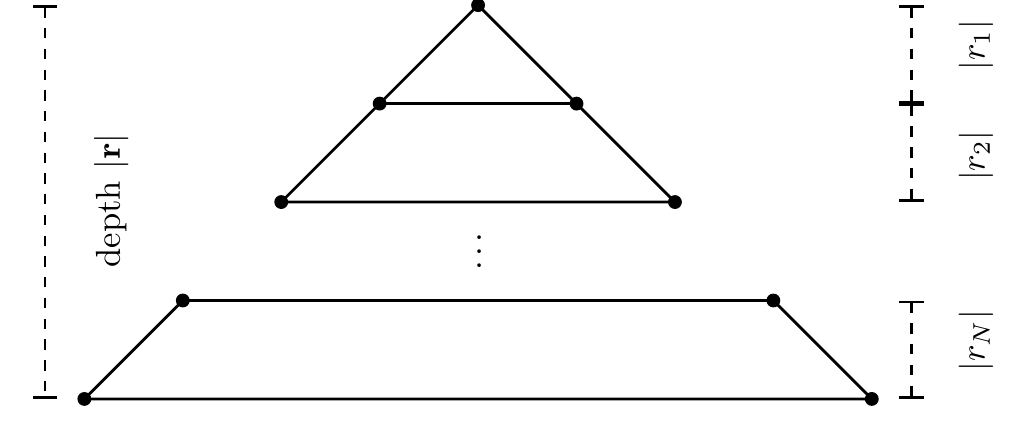}
\caption{Layerwise attaching to form $\T$.}
\label{figure:ErrRedBinaryTree2}
\vspace{-.2cm}
\end{wrapfigure}

% ----------------------
%\begin{Thm}\label{Theorem:ErrorGridReduction}
\begin{theorem}\label{Theorem:ErrorGridReduction}
Let $W^{\textnormal{pol}}(\cdot,t) \approx W(\cdot,t)$  
be the Chebyshev interpolation polynomial of the potential $W$ on $\Omega = [-L,L]^N$ over $\Rc(R)$ for fixed $L$. Let $\Kc(K)$ be 
a hyperbolically reduced index set with $K \gg R$. 
Then, under assumption \eqref{align:GeneralDecayCondition} on $v \in \C^{|\Kc|}$ (i.e., componentwise decay of order $\beta \in \N$), 
the error due to index set reduction behaves as
\begin{align*}
 \left|\left(\left(W^{\textnormal{pol}}(X_{\Kc},t)-\Wc_{\Kc,\textnormal{pol}}^{\textnormal{GH}(K)}(t)\right)v\right)_{\bfj}\right| \leq
           C(N,\Rc,W,\beta,L,t)K^{-\beta},\quad\quad\bfj \in \Kc.
\end{align*}
The matrix $W^{\textnormal{pol}}(X_{\Kc},t)$ results from formally inserting the hyperbolically reduced 
coordinate matrices $X_{\Kc}^{(l)}$  into the polynomial (see Section \ref{section:FastAlgorithm}) and $\Wc_{\Kc,\textnormal{pol}}^{\textnormal{GH}(K)}(t)$ is defined as in 
Section \ref{subsection:GaussHermiteQuadrature}.  The constant $C(N,\Rc,W,\beta,L,t)$ behaves as in \eqref{align:CRed}, see below, and depends only on $N, \Rc$, the regularity of $W$, $\beta$, $L$, and the time $t$.%\hfill\dollar
\end{theorem}
%\end{Thm}
\emph{Proof.} As in the previous theorem, we consider a partition of the error
\begin{align*}
 E^{\text{red}}(v)_\bfj = \sum_{\bfr \in \Rc} \alpha_{\bfr}(t) \left[\left(T_{\bfr}\left(\frac{1}{L}X_{\Kc}\right)v\right)_{\bfj}-\left(T_{\bfr}\left(\frac{1}{L}X_{\text{full}}\right)\Omega_+(v)\right)_{\bfj}\right], \quad \bfj \in \Kc.
\end{align*}
\textit{Construction of binary trees:} For fixed $l \in \{1,\dots,N\}$, $r_l \in \{1,\dots,R\}$ and $\bfj \in \Kc_{\text{full}}$, applying the Chebyshev recursion \eqref{align:ChebyshevRecurrence} together with \eqref{align:DirectOperationWithX} 
yields an expansion
\begin{align*}%\label{align:ErrorGridReductionConstructionBinaryTree}
 T_{r_l,\bfj}=\frac{\sqrt{2j_l}}{L}\sum_{\substack{s=1,\\s\text{ odd}}}^{r_l-1}(-1)^{s+1}T_{r_l-s, \bfj-\bfe_l}\!\!+\frac{\sqrt{2(j_l+1)}}{L}\sum_{\substack{s=1,\\s\text{ odd}}}^{r_l-1}(-1)^{s+1}T_{r_l-s, \bfj+\bfe_l} + \tau_{r_l,\bfj},
\end{align*}
where we use the abbrevations $T_{r_l,\bfj} = \left(T_{r_l}\left(\frac{1}{L}X_{\text{full}}^{(l)}\right)v\right)_{\bfj}$ and
\begin{align*}
 \tau_{r_l,\bfj} = \begin{cases}
               \pm v_{\bfj},& r_l \text{ even and } 4 \mid (\nmid)\,r_l,\\
               \pm \frac{1}{L}\left(\sqrt{\frac{j_l}{2}}v_{\bfj-\bfe_l} + \sqrt{\frac{j_l+1}{2}}v_{\bfj+\bfe_l}\right),& r_l \text{ odd and } 4 \mid (\nmid)\,(r_l+1).
              \end{cases}
\end{align*}
Repeated application allows for a binary tree expansion $\T^{(l)}(\bfj)$ of depth $r_l$ as given 
in Figure \ref{figure:ErrRedBinaryTree} (coefficients omitted). With each left or right descent, all indices of newly expanded $T$-terms have their $l$th component  
reduced or increased by $1$, respectively. Clearly, leaves carry sums of at most $\frac{1}{2}r_l(r_l-1)$ terms.
Starting from $\bfj \in \Kc_{\text{full}}$, a binary tree $\T_{\text{full}}(\bfj)$ for the product
\begin{align*}
 \left(T_{\bfr}\left(\frac{1}{L}X_{\text{full}}\right)v\right)_\bfj =  \left( T_{r_1}\left(\frac{1}{L}X_{\text{full}}^{(1)}\right) \cdot \left({\dots}\left(T_{r_N}\left(\frac{1}{L}X_{\text{full}}^{(N)}\right)v\right){\dots}\right)\right)_\bfj
\end{align*}
is then obtained by attaching to each term in each leaf of $\T^{(l)}$ analogously defined trees $\T^{(l+1)}$ starting from $l = 1$, such that 
a leaf in layer $l$ is a root of a subtree in layer $l+1$, see the pattern given in Figure \ref{figure:ErrRedBinaryTree2}, where the topmost and 
lowermost layers are numbered $1$ and $N$, respectively. Along the path to a proper leaf in layer $N$ (an $N$-\emph{leaf}), let $\lambda_l$ and $\rho_l$ denote the 
number of left and right descents in layer $l$, respectively. Starting from $\bfj$, in layer $l$, only the $l$th component of $\bfj$ is changed.\\

The same considerations apply with $X_{\Kc}^{(l)}$ in place of $X_{\text{full}}^{(l)}$, yielding an analogously defined binary tree $\T_{\Kc}(\bfj)$ starting from $\bfj \in \Kc$. 
We consider the difference tree $\D(\bfj) = \T_{\Kc}(\bfj)-\T_{\text{full}}(\bfj)$ for $\bfj \in \Kc$, using the vectors $v \in \C^{|\Kc|}$ and $\Omega_+(v)$. 
If an index does not belong to $\Kc_{\text{full}}$ or $\Kc$, we say that the corresponding term \emph{vanishes} in $\T_{\text{full}}(\bfj)$ or $\T_{\Kc}(\bfj)$, respectively. 
A term in the difference tree $\D(\bfj)$ \emph{vanishes} 
if corresponding terms in $\T_{\text{full}}(\bfj)$ and $\T_{\Kc}(\bfj)$ vanish or do not vanish both at the same time. We state the following obvious, yet important 
observations: Terms with an index belonging not even to $\Kc_{\text{full}}$ vanish  in $\D(\bfj)$ anyway. 
An $N$-leaf does not vanish in $\D(\bfj)$ 
if and only if, along the path connecting it to the root $v_\bfj$, there is at least one node belonging to $\Kc_{\text{full}} \setminus \Kc$. 
As in Section \ref{subsection:ErrorQuadrature}, we examine non-vanishing $N$-leaves in $\D(\bfj)$.\\
\textit{Characterization of non-vanishing leaves:} We consider a root index $\bfm$ in layer $l$, where $m_l = j_l$. We have the following requirements for a term depending on $\bfm$ not to vanish in $\T_{\text{full}}(\bfj)$ or $\T_{\Kc}(\bfj)$, respectively:
\begin{itemize}
 \item[($\T_{\text{full}}$)] For the $l$th component index, it is required that $0 \stackrel{!}{\leq} j_l + \rho_l - \lambda_l \leq j_l + r_l - 2\lambda_l \stackrel{!}{\leq} K$,
which gives the bounds
\begin{align*}
 \lambda_{\min}^l(\bfr,\bfj) = \left\lfloor \frac{j_l + r_l -K}{2} \right\rfloor \leq \lambda_l \leq \left\lfloor \frac{j_l + r_l}{2} \right\rfloor = \lambda_{\max}^l(\bfr,\bfj).
\end{align*}
\item[($\T_{\Kc}$)] The upper bound is the same as in $(\T_{\text{full}})$.  By the definition of $\Kc$, one needs
\begin{align}\label{align:HypCrossCondition}
  j_l + r_l - 2\lambda_l + 1\stackrel{!}{\leq} (K + 1) \Big(\prod_{\substack{i=1\\ i \neq l}}^N (1 + m_i)\Big)^{-1}.
\end{align}
From $\bfm \in \Kc$, it follows that
\begin{align*}
   1 + j_l \leq (K + 1) \Big(\prod_{\substack{i=1\\ i \neq l}}^N (1 + m_i)\Big)^{-1},
\end{align*}
thus, $r_l - 2\lambda_l \stackrel{!}{\leq} 0$, i.e., the leaf in $\T_{\Kc}(\bfj)$ does not vanish for
\begin{align}\label{align:HypCrossCondition2}
 \lambda_l \geq \left\lceil \frac{r_l}{2}\right\rceil = \lambda_{\min,\text{hyp}}^l(\bfr).
\end{align}
\end{itemize}
Non-vanishing leaves in $\D(\bfj)$ satisfy ($\T_{\text{full}}$), but not the more restrictive ($\T_{\Kc}$). The converse is not true, since a leaf violating $\T_{\Kc}$ might still 
fulfill \eqref{align:HypCrossCondition}, and thus vanish in $\D(\bfj)$. We consider the simpler 
condition \eqref{align:HypCrossCondition2}. Obviously, $\lambda_{\min}^l \leq \lambda_{\min,\text{hyp}}^l$.\\
\textit{Number of non-vanishing leaves:} Summing up as in Section \ref{subsection:ErrorQuadrature}, we have at most
\begin{align*}
 \sum_{s=\lambda_{\min}^l}^{\lambda_{\max}^l} \binom{r_l}{s} - \sum_{s=\lambda_{\min,\text{hyp}}^l}^{\lambda_{\max}^l} \binom{r_l}{s} = \sum_{s=\lambda_{\min}^l}^{\lambda_{\min,\text{hyp}}^l-1} \binom{r_l}{s} = a_l(\bfr,\bfj)
\end{align*}
non-vanishing $l$-leaves. We use \eqref{align:SumBinomialCoeff} with $a = \lambda_{\min}^l,\,b = \lambda_{\min,\text{hyp}}^l-1,\,c = r_l$. The fact $\lambda_{\min,\text{hyp}}^l = \left\lceil\frac{r_l}{2}\right\rceil$ 
implies $2(\lambda_{\min,\text{hyp}}^l-1) < r_l$, hence,
\begin{align*}
 a_l(\bfr,\bfj) < \binom{r_l}{\lfloor \frac{r_l}{2}-1\rfloor }\frac{r_l}{2}.
\end{align*}
\textit{Error accumulation:} Along the path to any $l$-leaf, the most unfavorable weight
\begin{align*}
 b_l(\bfr,\bfj) = 2^{r_l/2} \prod_{s=1}^{r_l} \frac{(j_l+s)^{1/2}}{L}% = \Oc\left(2^{r_l/2} S^{r_l/2}\right),
\end{align*}
 comes from descending right only. Using \eqref{align:GeneralDecayCondition}, the largest $N$-leaf is bounded by
\begin{align*}
 c(\bfr,\bfj) \leq \prod_{l=1}^N \frac{1}{2}r_l(r_l-1) \left\{\max_{-r_l \leq s \leq r_l} \left|v_{\bfj - s\bfe_l}\right| \right\} = \Oc\left(\prod_{l=1}^N \frac{1}{2}r_l(r_l-1) (j_l-r_l)^{-\beta}\right).
\end{align*}
A path in $\D(\bfj)$ does not vanish only in case $\bfj + \bfr \in \Kc_{\text{full}} \setminus \Kc$, thus,
\begin{align*}
 K\!+\!2\!\leq\!\prod_{l=1}^N (1\!+\!j_l\!+\!r_l)\!=\!\prod_{l=1}^N (j_l\!-\!r_l) \cdot \prod_{l=1}^N \left(1 + \frac{1\!+\!2r_l}{j_l-r_l}\right)\!\leq\!\prod_{l=1}^N (j_l\!-\!r_l) \cdot 2^N \prod_{l=1}^N (1\!+\!r_l).
\end{align*}
Hence, the error over all layers is bounded as
\begin{align*}%\label{align:ErrorGridReductionGeneralBound}
 |E_\bfj^\bfr| &\leq  \prod_{l=1}^N \left\{a_l(\bfr,\bfj) \cdot b_l(\bfr,\bfj) \right\} \cdot c(\bfr,\bfj)\\ %= \Oc\left(\prod_{l=1}^N a_l(\bfr,\bfj) \cdot b_l(\bfr,\bfj)\cdot r_l (j_l-r_l)^{-\beta}\right)\\
 %&= \Oc\left(2^{\beta N}\prod_{l=1}^N \left\{ 	\binom{r_l}{\lfloor\frac{r_l}{2}\!-\!1\rfloor}\frac{r_l}{2} \cdot 2^{r_l/2}\prod_{s=1}^{r_l} \frac{(j_l\!+\!s)^{1/2}}{L} \cdot \frac{1}{2}r_l(r_l\!-\!1)(1\!+\!r_l)^{\beta}(K\!+\!2)^{-\beta}	\right\}\right)\\
 &= \Oc\left(2^{(\beta-2) N}\prod_{l=1}^N \left\{ 	\binom{r_l}{\lfloor\frac{r_l}{2}\!-\!1\rfloor} \cdot 2^{r_l/2}\prod_{s=1}^{r_l} \frac{(j_l\!+\!s)^{1/2}}{L} \cdot r_l^{\beta+3}(K\!+\!2)^{-\beta}	\right\}\right).
\end{align*}
Finally, we sum up and set
\begin{align}\label{align:CRed}
 \hspace{-.4cm}C(N,\Rc,W,\beta,L,t)\!=\!\!\sum_{\bfr \in \Rc}\!\!\alpha_{\bfr}(t)2^{(\beta\!-\!2) N}\prod_{l=1}^N \left\{ 	\binom{r_l}{\lfloor\frac{r_l}{2}\!-\!1\rfloor} 2^{r_l/2}\prod_{s=1}^{r_l} \frac{(j_l\!+\!s)^{1/2}}{L} r_l^{\beta+3}\right\}.
\end{align}\hfill\dollar
\textit{Remarks:} According to the choice of $\bfj \in \Kc$ or $\Rc(R)$, the above error estimate might improve. If there is more than one large component in $\bfj$, say $n(\Rc,\bfj)= \min_{\bfr \in \Rc} n(\bfr,\bfj) \geq 2$, where $n(\bfr,\bfj) \in \{1,\dots,N\}$ is the number of components $j_l$ such that 
 $K \approx j_l \gg r_l$, we get $c(\bfr,\bfj) = \Oc(K^{-n(\Rc,\bfj)\beta})$. 
 On the other hand, if $\bfj + \bfr \in \Kc$ (i.e., only small index components), all paths in $\D(\bfj)$ cancel out and the error $E^{\bfr}_{\bfj}$ vanishes.
%%%%%%%%%%%%%%%%%%%%%%%%%%%%%%%%%%%%%%%%%%%%%%%%%%%%%%%%%%%%%%%%%%%%%%%%%%%%%%%%%%%%%%
%%%%%%%%%%%%%%%%%%%%%%%%%%%%%%%%%%%%%%%%%%%%%%%%%%%%%%%%%%%%%%%%%%%%%%%%%%%%%%%%%%%%%%
\section{Numerical experiments}\label{section:NumericalExperiments}
All figures have been obtained on a desktop computer with an Intel Core 2 
Duo E8400 3.00 GHz processor with 4 GB RAM.
%%%%%%%%%%%%%%%%%%%%%%%%%%%%%%%%%%%%%%%%%%%%%%%%%%%%%%%%%%%%%%%%%%%%%%%%%%%%%%%%%%%%%%
\subsection{Local errors due to quadrature and index set reduction}
\begin{figure}[!ht]
\centering
\begin{minipage}{.6\linewidth}\hspace{-.3cm}
\includegraphics[scale=0.25]{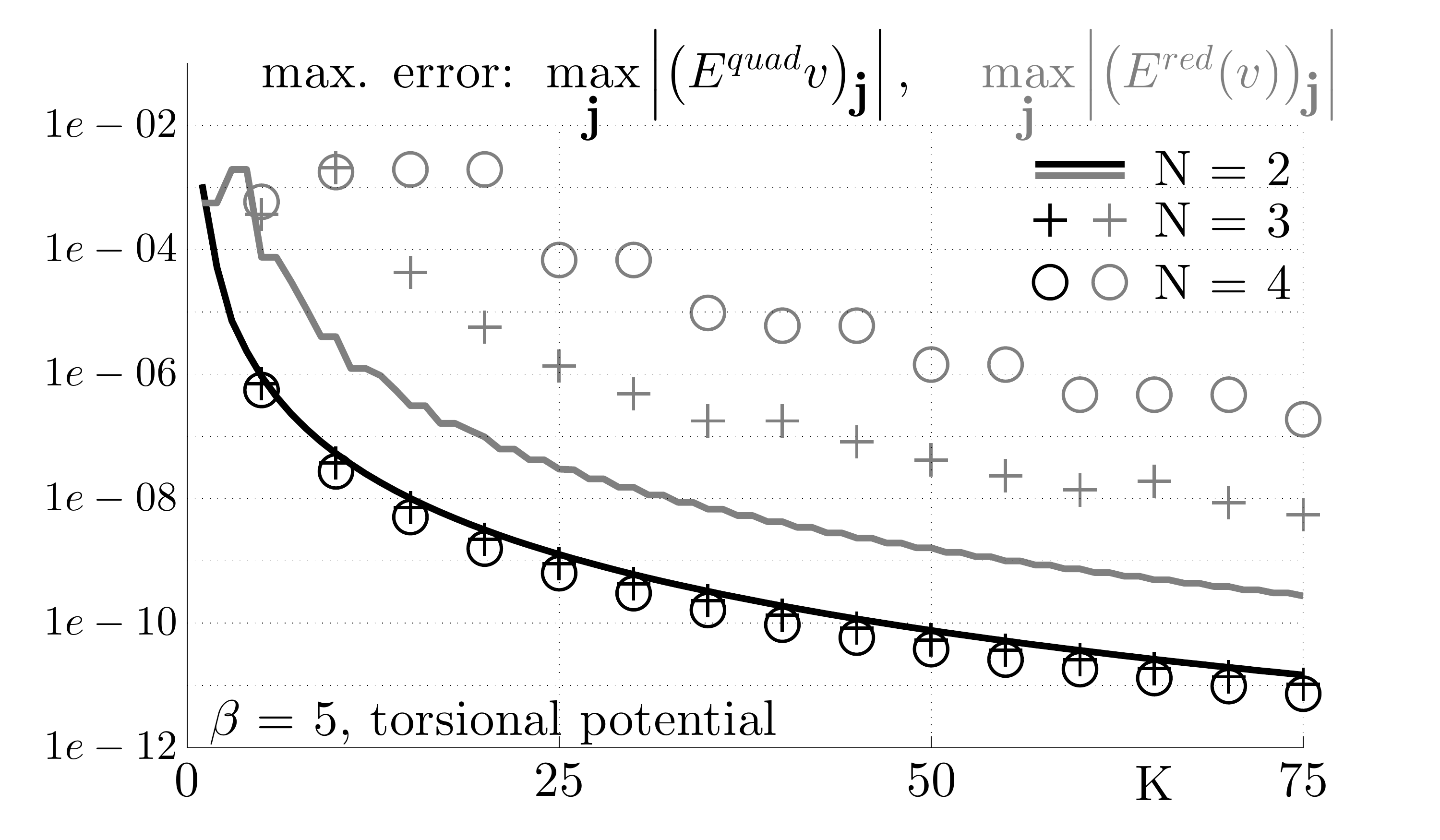}
\end{minipage}
\begin{minipage}{.39\linewidth}\hspace{-1.0cm}
 {%\tabulinesep=1.2mm
 \begin{tabular}{|c|c|c|c|}
  \hline \small$K \downarrow$& \small$N = 2$& \small$N=3$& \small$N=4$\\
  \hline \hline \small 25& \scriptsize $1.270\mh09$& \scriptsize $8.978\mh10$& \scriptsize $6.347\mh10$\\
         \hline \small 50& \scriptsize $7.629\mh11$& \scriptsize $5.392\mh11$& \scriptsize $3.812\mh11$\\
         \hline \small 75& \scriptsize $1.466\mh11$& \scriptsize $1.037\mh11$& \scriptsize $7.328\mh12$\\
\hline \hline \multirow{1}{*}{\scriptsize time}& \scriptsize 19.0 secs& \scriptsize 5.9 min& \scriptsize 55.8 min\\
\hline \hline \small \textcolor{gray}{25}& \scriptsize \textcolor{gray}{$2.975\mh08$}& \scriptsize \textcolor{gray}{$1.344\mh06$}& \scriptsize \textcolor{gray}{$6.834\mh05$}\\
\hline		\small \textcolor{gray}{50}& \scriptsize \textcolor{gray}{$1.621\mh09$}& \scriptsize \textcolor{gray}{$4.207\mh08$}& \scriptsize \textcolor{gray}{$1.426\mh06$}\\
\hline 		\small \textcolor{gray}{75}& \scriptsize \textcolor{gray}{$2.729\mh10$}& \scriptsize \textcolor{gray}{$5.540\mh09$}& \scriptsize \textcolor{gray}{$1.879\mh07$}\\
\hline \hline \multirow{1}{*}{\scriptsize \textcolor{gray}{time}}& \scriptsize \textcolor{gray}{0.002 secs}& \scriptsize \textcolor{gray}{0.011 secs}& \scriptsize \textcolor{gray}{0.047 secs}\\
\hline
 \end{tabular}}
\end{minipage}
\caption{Errors $E^{\text{quad}}v$ (black) and $E^{\text{red}}(v)$ (\textcolor{gray}{gray}) for the torsional potential \eqref{align:TorsionalPotential} 
($L\!=\!16$, $S\!=\!1$, $R\!=\!8$, $\beta\!=\!5$). The solid lines represent $N\!=\!2$, selected errors in cases $N\!=\!3, 4$ 
are indicated by plus signs 
and circles, respectively. Increasing $N$ worsens the factor $C(N,\Rc,W,\beta,L)$ in $E^{\text{red}}$, see \eqref{align:CRed}. 
In the rows labeled ``time'', computation times for $\mathcal{W}_{\Kc,\text{pol}}^{\text{GH}(K)}v$ (assembly and multiplication) and for the fast algorithm in case $K\!=\!75$ are shown. 
As for $N\!=\!4$, assembling the matrix (plus operating on a vector) takes almost one hour -- even with a reduced index set.}
\label{figure:LokalerFehlerN2bis5}
\vspace{-.2cm}
\end{figure}

\begin{figure}[!ht]
\floatbox[{\capbeside\thisfloatsetup{capbesideposition={right,top},capbesidewidth=5.8cm}}]{figure}[\FBwidth]
{\caption{Errors $\left(E^{\text{red}}(v)\right)_\bfj$ due to index set reduction for a torsional potential with $N\!=\!2$ and  different choices of $K$ (as above, $L\!=\!16$, $S\!=\!1$, $R\!=\!8$, $\beta\!=\!3$). Each entry represents an error vector component. Errors being small with respect to the largest observed error component $e_{\max} \approx 5.264\mh05$ are simply indicated by a dot,  
indices carrying larger errors are indicated by a grey box. The darker the box, the closer the error to $e_{\max}$. The pictures corresponding 
to $K\!=\!30, 40$ show an enlarged view. The symbol \# points to the number of large error components. The errors 
decrease with growing $K$ as indicated by increasingly lighter boxes and are concentrated in the region with only ``intermediate'' index components.}\label{figure:FehlervektorSchnellerAlgorithmusN2}}
{\includegraphics[trim=1.3cm 1cm 1.4cm .4cm,clip,scale=0.300]{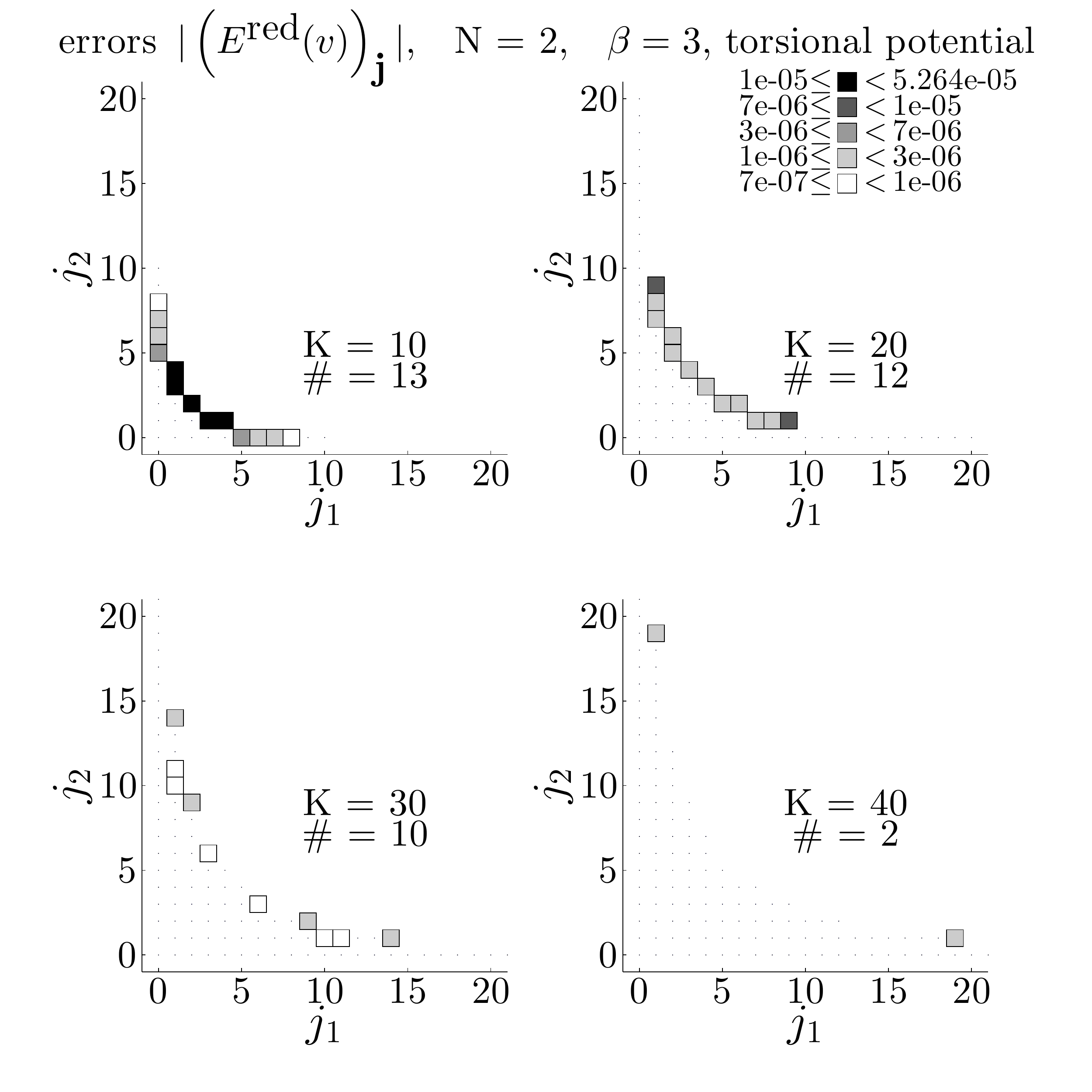}}
\end{figure}

Let $\Kc = \Kc(K)$ be a hyperbolically reduced $N$-dimensional index set. We illustrate the errors
\begin{align*}
 E^{\text{quad}}v = \left(\Wc_{\Kc,\text{pol}} - \Wc_{\Kc,\text{pol}}^{\text{GH}(K)}\right)v,\quad\quad  E^{\text{red}}(v) =  \left(W^{\text{pol}}(X_{\Kc}) - \Wc_{\Kc,\text{pol}}^{\text{GH}(K)}\right)v
\end{align*}
due to quadrature  and index set reduction as given in Theorems \ref{Theorem:ErrorQuadratureMatrixVector} and \ref{Theorem:ErrorGridReduction}, respectively, 
for different choices of $N$ and $K$, see Figure \ref{figure:LokalerFehlerN2bis5}. 
In both cases, the chosen potential is the aforementioned stretched torsional potential as 
given in \eqref{align:TorsionalPotential}, i.e.,  $W(x) = \sum_{l=1}^N \left(1 - \cos(x_l/L)\right)$ with $L=16$ approximated by Chebyshev interpolation with $R = 8$ nodes on each axis as in Section \ref{subsection:Complexity}. 
For the vector $v \in \C^{|\Kc|}$ to exhibit a decay behavior according to \eqref{align:GeneralDecayCondition}, we set
\begin{align*}
 v_\bfk = \prod_{l=1}^N\max(k_l,1)^{-\beta},\qquad \beta = 5, 
\end{align*}
and then normalize such that $\|v\|_2 = 1$. As explained in Section \ref{subsection:ErrorGridReduction}, for $K$ 
being sufficiently large, the error $\left(E^{\text{red}}(v)\right)_\bfj$ 
decreases the faster the more components $j_l$ of $\bfj$ are large with respect to $R$, 
see Theorem \ref{Theorem:ErrorGridReduction} and the remarks thereafter. 
Figure \ref{figure:FehlervektorSchnellerAlgorithmusN2} illustrates this decay behavior in the individual components of $E^{\text{red}}(v)$ for $N = 2$ and $\beta = 3$. The matrix $\Wc_{\Kc,\text{pol}}^{\text{GH}(K)}$ is assembled as explained in Section \ref{subsection:GaussHermiteQuadrature}. 
For a full index set $\Kc_{\text{full}}$, Lemma \ref{Lemma:QuadratureEqualsInsertion} has been confirmed  numerically 
for all our choices of $K$ and $N$.
%-------------------------------
%---------------------------
%----------------------------
%%%%%%%%%%%%%%%%%%%%%%%%%%%%%%%%%%%%%%%%%%%%%%%%%%%%%%%%%%%%%%%%%%%%%%%%%%%%%%%%%%%%%%
\subsection{Perturbed Lanczos process}\label{subsection:perturbedLanczos}
\begin{figure}[!ht]
\centering
\begin{minipage}{.47\linewidth}\hspace{-.4cm}
\includegraphics[trim=.0cm .6cm .0cm .0cm,clip,scale=0.255]{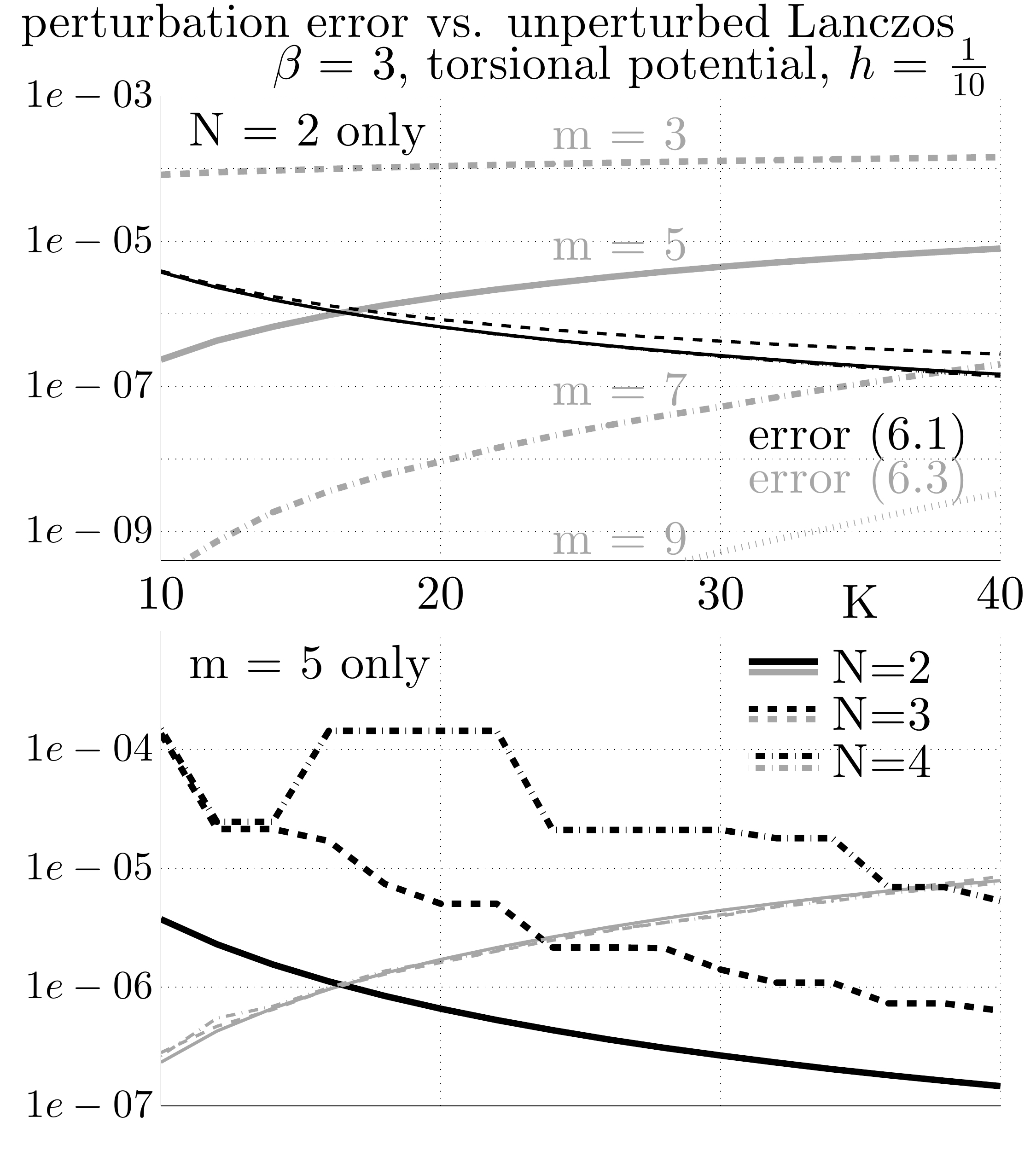}
\end{minipage}
\begin{minipage}{.52\linewidth}\hspace{-0.8cm}
\vspace{0.4cm}
 %\tabulinesep=1.2mm
 \begin{tabular}{|c|c|c|c|c|}
  \hline $K\!\rightarrow$				 & $10$						 & $20$						 & $30$						 & $40$\\
  \hline \hline \scriptsize \textcolor{gray}{$N\!=\!2$}& \scriptsize \textcolor{gray}{$8.196 \mh 05$}& \scriptsize \textcolor{gray}{$1.077 \mh 04$}& \scriptsize \textcolor{gray}{$1.270 \mh 04$}& \scriptsize \textcolor{gray}{$1.427 \mh 04$}\\
		\scriptsize \textcolor{color4}{$m\!=\!3$}& \scriptsize \textcolor{color4}{$3.878 \mh 06$}& \scriptsize \textcolor{color4}{$8.279 \mh 07$}& \scriptsize \textcolor{color4}{$4.181 \mh 07$}& \scriptsize \textcolor{color4}{$2.779 \mh 07$}\\
         \hline \scriptsize \textcolor{gray}{$N\!=\!2$}& \scriptsize \textcolor{gray}{$2.324 \mh 07$}& \scriptsize \textcolor{gray}{$1.709 \mh 06$}& \scriptsize \textcolor{gray}{$4.415 \mh 06$}& \scriptsize \textcolor{gray}{$7.881 \mh 06$}\\
		\scriptsize \textcolor{color4}{$m\!=\!5$}& \scriptsize \textcolor{color4}{$3.731 \mh 06$}& \scriptsize \textcolor{color4}{$6.565 \mh 07$}& \scriptsize \textcolor{color4}{$2.650 \mh 07$}& \scriptsize \textcolor{color4}{$1.464 \mh 07$}\\
         \hline \scriptsize \textcolor{gray}{$N\!=\!2$}& \scriptsize \textcolor{gray}{$2.522 \mh 10$}& \scriptsize \textcolor{gray}{$9.192 \mh 09$}& \scriptsize \textcolor{gray}{$5.211 \mh 08$}& \scriptsize \textcolor{gray}{$2.003 \mh 07$}\\
		\scriptsize \textcolor{color4}{$m\!=\!7$}& \scriptsize \textcolor{color4}{$3.730 \mh 06$}& \scriptsize \textcolor{color4}{$6.533 \mh 07$}& \scriptsize \textcolor{color4}{$2.588 \mh 07$}& \scriptsize \textcolor{color4}{$1.374 \mh 07$}\\
	 \hline \scriptsize \textcolor{gray}{$N\!=\!2$}& \scriptsize \textcolor{gray}{$1.179 \mh 13$}& \scriptsize \textcolor{gray}{$2.973 \mh 11$}& \scriptsize \textcolor{gray}{$5.241 \mh 10$}& \scriptsize \textcolor{gray}{$3.354 \mh 09$}\\
		\scriptsize \textcolor{color4}{$m\!=\!9$}& \scriptsize \textcolor{color4}{$3.730 \mh 06$}& \scriptsize \textcolor{color4}{$6.533 \mh 07$}& \scriptsize \textcolor{color4}{$2.587 \mh 07$}& \scriptsize \textcolor{color4}{$1.371 \mh 07$}\\
\hline	\hline	\scriptsize \textcolor{gray}{$N\!=\!3$}& \scriptsize \textcolor{gray}{$2.796 \mh 07$}& \scriptsize \textcolor{gray}{$1.612 \mh 06$}& \scriptsize \textcolor{gray}{$4.062 \mh 06$}& \scriptsize \textcolor{gray}{$8.531 \mh 06$}\\
		\scriptsize \textcolor{color4}{$m\!=\!5$}& \scriptsize \textcolor{color4}{$1.363 \mh 04$}& \scriptsize \textcolor{color4}{$5.035 \mh 06$}& \scriptsize \textcolor{color4}{$1.405 \mh 06$}& \scriptsize \textcolor{color4}{$6.310 \mh 07$}\\
\hline 		\scriptsize \textcolor{gray}{$N\!=\!4$}& \scriptsize \textcolor{gray}{$2.628 \mh 07$}& \scriptsize \textcolor{gray}{$1.666 \mh 06$}& \scriptsize \textcolor{gray}{$3.966 \mh 06$}& \scriptsize \textcolor{gray}{$7.500 \mh 06$}\\
		\scriptsize \textcolor{color4}{$m\!=\!5$}& \scriptsize \textcolor{color4}{$1.493 \mh 04$}& \scriptsize \textcolor{color4}{$1.437 \mh 04$}& \scriptsize \textcolor{color4}{$2.103 \mh 05$}& \scriptsize \textcolor{color4}{$5.313 \mh 06$}\\
\hline
 \end{tabular}
\end{minipage}
\caption{Errors \eqref{align:ErrorPerturbedLanczosAnalysis} (\textcolor{color4}{black}) and \eqref{align:ErrorUnperturbedLanczos} (\textcolor{gray}{gray}) as functions of $K$ for fixed $h\!=\!1/10$ and various choices of $N$ and $m$ (torsional potential, $L\!=\!16$, $S\!=\!1$, $R\!=\!8$, $\beta\!=\!3$). 
Upper figure: $N\!=\!2$, $m\!=\!3,\!\dots\!,9$ (dashed, solid, chain dotted, and dotted line, respectively). Lower figure: $N\!=\!2,3,4$, $m\!=\!5$ (solid, dashed, and chain dotted line, respectively). 
In each cell of the table, the upper and lower figure corresponds to \eqref{align:ErrorUnperturbedLanczos} and \eqref{align:ErrorPerturbedLanczosAnalysis}, respectively.}
\label{figure:perturbation_lanczos_classvspert}
\end{figure}

For a problem of the form $i\dot{y}(t) = Ay(t)$, we approximate the matrix exponential $\exp(-ihA)v$ 
using an $m$-step Lanczos process. 
% 
% We approximate the  . Again, $W$ is the 
% above torsional potential ($L = 16$, $S = 1$, Chebyshev interpolation, $R = 8$) and $v$ is given as above. 
In each Lanczos step, applying the fast algorithm 
\begin{align*}
 (Av_k)^{fast} = Av_k - \left(Av_k - (Av_k)^{fast}\right) = Av_k - f_k
\end{align*}
instead of $Av_k$ produces perturbed basis vectors and coefficients $\tilde{V}_m$ and $\tilde{T}_m$, respectively. This yields
 $A = \tilde{V}_m\tilde{T}_m\tilde{V}_m^* + F_m\tilde{V}_m^*$, 
where $F_m = \left(f_1|\dots|f_m\right)$ are the perturbations. Thus, the relation to the unperturbed counterparts $V_m$ and $T_m$ reads $V_mT_mV_m^* = \tilde{V}_m\tilde{T}_m\tilde{V}_m^* + F_m\tilde{V}_m^*$. 
We approximate $\exp(-ihA)v \approx \tilde{V}_m \exp(-ih\tilde{T}_m)e_1$, 
and,  by a sensitivity analysis for the matrix exponential as done in \cite{VL77}, the local error is readily seen to be
\begin{align}\label{align:ErrorPerturbedLanczosAnalysis}
 \left\|V_m \exp(-ihT_m)e_1 - \tilde{V}_m \exp(-ih\tilde{T}_m)e_1\right\|_2 \leq h \| F_m\|_2 \exp(h(\|A\|_2 + \|F_m\|_2)).
\end{align}
We consider $A = \mathcal{D}_{\Kc} + \Wc_{\Kc,\text{pol}}$, see \eqref{align:ODECoefficients}, where $\mathcal{D}_{\Kc}$ is the diagonal matrix given in Section \ref{subsection:GalerkinAnsatz}, 
and the underlying potential is again the torsional potential,
\begin{align}\label{align:TorsionalPotentialwithHO}
W(x) = \sum_{l=1}^N \left(1 - \cos(x_l/L)\right) - \frac{1}{2}\sum_{l=1}^N x_l^2, \quad x \in \Omega,
\end{align}
 (Chebyshev interpolation as in Section \ref{subsection:Complexity}, $L=16$, $R=8$). The linear decay of the error \eqref{align:ErrorPerturbedLanczosAnalysis} with respect to $h$ is shown in Figure \ref{figure:perturbation_lanczos_errorvsh}. 
Additionally, for a fixed choice of $h$, the error is seen to become worse for constant $K$ and increasing $N$ (as predicted by Theorem \ref{Theorem:ErrorGridReduction}),
and it becomes arbitrarily small for constant $N$ and a sufficiently large choice of $K$. 
We apply only 5 Lanczos steps in each time step.
%-----------------------
\begin{figure}[!ht]
\floatbox[{\capbeside\thisfloatsetup{capbesideposition={right,top},capbesidewidth=4cm}}]{figure}[\FBwidth]
{\begin{tabular}{|c|c|c|c|c|}
  \hline $h\!\rightarrow$& $1/10$& $1/20$& $1/40$& $1/80$\\
  \hline \hline \scriptsize $N\!=\!2, K\!=\!10$& \scriptsize $3.731 \mh 06$& \scriptsize $1.874 \mh 06$& \scriptsize $9.380 \mh 07$& \scriptsize $4.691 \mh 07$\\
         \hline \scriptsize $N\!=\!2, K\!=\!20$& \scriptsize $6.565 \mh 07$& \scriptsize $3.285 \mh 07$& \scriptsize $1.644 \mh 07$& \scriptsize $8.222 \mh 08$\\
         \hline \scriptsize $N\!=\!2, K\!=\!30$& \scriptsize $2.650 \mh 07$& \scriptsize $1.305 \mh 07$& \scriptsize $6.526 \mh 08$& \scriptsize $3.264 \mh 08$\\
\hline 		\scriptsize $N\!=\!2, K\!=\!40$& \scriptsize $1.464 \mh 07$& \scriptsize $6.956 \mh 08$& \scriptsize $3.468 \mh 08$& \scriptsize $1.734 \mh 08$\\
\hline	\hline	\scriptsize $N\!=\!3, K\!=\!40$& \scriptsize $6.310 \mh 07$& \scriptsize $3.151 \mh 07$& \scriptsize $1.576 \mh 07$& \scriptsize $7.880 \mh 08$\\
\hline 	\hline	\scriptsize $N\!=\!4, K\!=\!40$& \scriptsize $5.313 \mh 06$& \scriptsize $2.661 \mh 06$& \scriptsize $1.331 \mh 06$& \scriptsize $6.658 \mh 07$\\
\hline
 \end{tabular}}
{\caption{Perturbation error \eqref{align:ErrorPerturbedLanczosAnalysis} depending linearly on $h$ for fixed $m\!=\!5$ and various choices of $N$ and $K$ 
(torsional potential, $L\!=\!16$, $S\!=\!1$, $R\!=\!8$, $\beta\!=\!3$).}
\label{figure:perturbation_lanczos_errorvsh}}
\end{figure}
%-----------------------
Note that if $m$ is chosen too large,  
the vectors $v_k$, $k \geq 2$, might fail to decay sufficiently fast, and the perturbation error might dominate the error due to Lanczos itself. 
A comparison of \eqref{align:ErrorPerturbedLanczosAnalysis} to the unperturbed Lanczos error
\begin{align}\label{align:ErrorUnperturbedLanczos}
 \left\|V_m \exp(-ihT_m)e_1 - \exp(-ihA)v\right\|_2
\end{align}
(see, e.g., \cite{L08}, Thm. III.2.10) is given in Figure \ref{figure:perturbation_lanczos_classvspert}, 
where we illustrate the error behavior for various choices of $N$, $K$, and $m$ and for a fixed choice of $h$. As the figures reveal, there is an antagonism: 
On the one hand, increasing $m$ improves \eqref{align:ErrorUnperturbedLanczos}, 
but \eqref{align:ErrorPerturbedLanczosAnalysis} might dominate unless $K$ is chosen sufficiently large. 
On the other hand, increasing $K$ requires a larger choice of $m$ for a decent Lanczos approximation (for fixed $h$). 
In our examples, moderate $m$ (say, $m =$ 5) is a good choice. Increasing $m$ is always possible, but yields additional costs without further reducing the overall error due to \eqref{align:ErrorPerturbedLanczosAnalysis} (unless $h$ is chosen smaller). 
Presently, we lack further analytical insight into how to relate $m$, $K$, and $h$ in an optimal way.
%-----------------------
%-----------------------
%%%%%%%%%%%%%%%%%%%%%%%%%%%%%%%%%%%%%%%%%%%%%%%%%%%%%%%%%%%%%%%%%%%%%%%%%%%%%%%%%%%%%%
\subsection{Time integration}
\label{subsection:timeintegration}
\begin{figure}[!ht]
\centering
\begin{minipage}{.99\linewidth}
\centering
 \includegraphics[scale=0.29,trim = 2cm .0cm 2cm .0cm,clip]{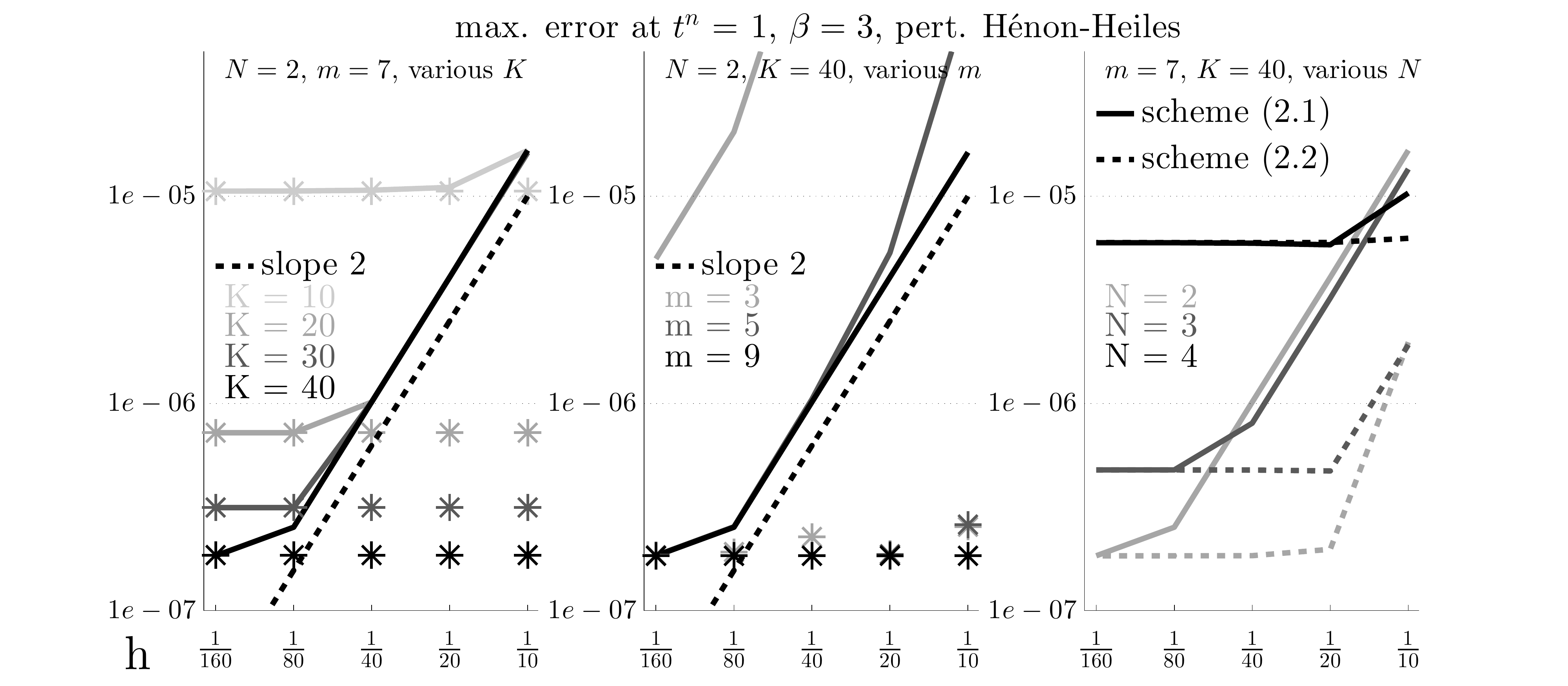}
\end{minipage}
\begin{minipage}{.99\linewidth}
\centering
 \begin{tabular}{|c|c|c|c|c|c|c|c|c|}
  \hline $N$& $m$& $K$&	\scriptsize scheme						       & $h\!=\!1/10$		   & $1/20$		       & $1/40$		    	   & $1/80$		       & $1/160$\\
  \hline \hline \scriptsize 2&\scriptsize 7&\scriptsize 10& 	\scriptsize\eqref{align:expMPR}& \scriptsize $1.665 \mh 05$& \scriptsize $1.103 \mh 05$& \scriptsize $1.069 \mh 05$& \scriptsize $1.060 \mh 05$& \scriptsize $1.058 \mh 05$\\
         \hline \scriptsize 2&\scriptsize 7&\scriptsize 20&	\scriptsize\eqref{align:expMPR}& \scriptsize $1.627 \mh 05$& \scriptsize $4.063 \mh 06$& \scriptsize $1.014 \mh 06$& \scriptsize $7.214 \mh 07$& \scriptsize $7.213 \mh 07$\\
         \hline \scriptsize 2&\scriptsize 7&\scriptsize 30&	\scriptsize\eqref{align:expMPR}& \scriptsize $1.623 \mh 05$& \scriptsize $4.062 \mh 06$& \scriptsize $1.014 \mh 06$& \scriptsize $3.141 \mh 07$& \scriptsize $3.141 \mh 07$\\
         \hline \scriptsize 2&\scriptsize 7&\scriptsize 40&	\scriptsize\eqref{align:expMPR}& \scriptsize $1.666 \mh 05$& \scriptsize $4.062 \mh 06$& \scriptsize $1.014 \mh 06$& \scriptsize $2.523 \mh 07$& \scriptsize $1.837 \mh 07$\\
\hline \hline	\scriptsize 2&\scriptsize 3&\scriptsize 40& 	\scriptsize\eqref{align:expMPR}& \scriptsize $1.568 \mh 02$& \scriptsize $4.047 \mh 03$& \scriptsize $2.763 \mh 04$& \scriptsize $2.044 \mh 05$& \scriptsize $4.987 \mh 06$\\
\hline		\scriptsize 2&\scriptsize 5&\scriptsize 40& 	\scriptsize\eqref{align:expMPR}& \scriptsize $9.066 \mh 05$& \scriptsize $5.330 \mh 06$& \scriptsize $1.045 \mh 06$& \scriptsize $2.537 \mh 07$& \scriptsize $1.843 \mh 07$\\
\hline 		\scriptsize 2&\scriptsize 9&\scriptsize 40& 	\scriptsize\eqref{align:expMPR}& \scriptsize $1.629 \mh 05$& \scriptsize $4.063 \mh 06$& \scriptsize $1.014 \mh 06$& \scriptsize $2.523 \mh 07$& \scriptsize $1.837 \mh 07$\\
\hline \hline 	\scriptsize 3&\scriptsize 7&\scriptsize 40& 	\scriptsize\eqref{align:expMPR}& \scriptsize $1.353 \mh 05$& \scriptsize $3.224 \mh 06$& \scriptsize $8.026 \mh 07$& \scriptsize $4.770 \mh 07$& \scriptsize $4.770 \mh 07$\\
\hline 		\scriptsize 4&\scriptsize 7&\scriptsize 40& 	\scriptsize\eqref{align:expMPR}& \scriptsize $1.036 \mh 05$& \scriptsize $5.830 \mh 06$& \scriptsize $5.933 \mh 06$& \scriptsize $5.960 \mh 06$& \scriptsize $5.967 \mh 06$\\
\hline 		\scriptsize 2&\scriptsize 7&\scriptsize 40& 	\scriptsize\eqref{align:2stageMagnus}& \scriptsize $1.973 \mh 06$& \scriptsize $1.974 \mh 07$& \scriptsize $1.840 \mh 07$& \scriptsize $1.837 \mh 07$& \scriptsize $1.837 \mh 07$\\
\hline 		\scriptsize 3&\scriptsize 7&\scriptsize 40& 	\scriptsize\eqref{align:2stageMagnus}& \scriptsize $1.918 \mh 06$& \scriptsize $4.719 \mh 07$& \scriptsize $4.769 \mh 07$& \scriptsize $4.770 \mh 07$& \scriptsize $4.770 \mh 07$\\
\hline 		\scriptsize 4&\scriptsize 7&\scriptsize 40& 	\scriptsize\eqref{align:2stageMagnus}& \scriptsize $6.266 \mh 06$& \scriptsize $5.975 \mh 06$& \scriptsize $5.969 \mh 06$& \scriptsize $5.969 \mh 06$& \scriptsize $5.969 \mh 06$\\
%\hline 		\scriptsize 4&\scriptsize 7&\scriptsize 40& 	\scriptsize\eqref{align:2stageMagnus}& \scriptsize $6.266 \mh 06$& \scriptsize $5.975 \mh 06$& \scriptsize $5.969 \mh 06$& \scriptsize $5.969 \mh 06$& \scriptsize $5.969 \mh 06$\\
\hline
 \end{tabular}
\end{minipage}
\caption{Propagation of \eqref{align:PropagatedEquation} with $W(x)$ as in \eqref{align:HenonHeilesExperiments} 
(Chebyshev interpolation, $L\!=\!16$, $S\!=\!1$, $R\!=\!3$, $\beta\!=\!3$) 
using the schemes \eqref{align:expMPR} and \eqref{align:2stageMagnus} as given in the table containing the observed errors \eqref{align:timePropagationError}. 
Corresponding errors \eqref{align:timePropagationPerturbationError} due to a perturbation of Lanczos as observed in the last time step are indicated by asterisks (not included in the table).}
\label{figure:ErrorTimeIntegration}
\vspace{-.2cm}
 \end{figure}

%We consider an instance of \eqref{align:SchroedingerEquation} where \eqref{align:ODECoefficientsPol} is a good approximation due to the solution being confined to a %given cube over integration time. 
We propagate \eqref{align:ODECoefficientsPol}, i.e., 
\begin{align}\label{align:PropagatedEquation}
\begin{aligned}
i \dot{c}_{\text{pol}}(t) &= \mathcal{D}_{\Kc} c_{\text{pol}}(t) + \mathcal{W}_{\Kc,\text{pol}}(t) c_{\text{pol}}(t),\\
\left(\tilde{c}_{\text{pol}}(0)\right)_{\bfk} &= \prod_{k_l\neq 0} k_l^{-\beta}, \quad \left(c_{\text{pol}}(0)\right)_{\bfk} = \left(\tilde{c}_{\text{pol}}(0)\right)_{\bfk} / \left\|\tilde{c}_{\text{pol}}(0)\right\|_2,
\end{aligned}
\end{align}
over $[0,1]$, where $\mathcal{D}_{\Kc}$ is the above diagonal matrix and the underlying potential $W$ is a stretched H\'{e}non-Heiles potential 
with a linear time-dependent perturbation,
\begin{align}\label{align:HenonHeilesExperiments}
 W(x,t) =\!\sum_{l=1}^{N-1}\!\left[(x_l/L)^2(x_{l+1}/L)\!-\!\frac{1}{3}(x_{l+1}/L)^3\right] -\!\sin^2(t)x_1 -\!\frac{1}{2}\sum_{l=1}^N\!x_l^2, \quad x \in \Omega,
\end{align}
(fully-indexed Chebyshev interpolation, $L = 16$, $R = 3$). This models the interaction of an atom or a molecule with a high-intensity CW laser in $x_1$-direction, 
see \cite{PKM94} (with a quantum harmonic oscillator in place of a HH-potential). 
%%%
Our aim is to show numerically the expected order of convergence with respect to $h$ and to illustrate the error behavior if $K$ and $m$ vary individually. 
%%%
Figure \ref{figure:ErrorTimeIntegration} shows the error
\begin{align}\label{align:timePropagationError}
 \max_{\bfj\in\Kc} \left|\left(c_{\text{pol}}^n-c_{\text{pol}}(t^n)\right)_{\bfj}\right|
\end{align}
at time $t^n = 1$ when using the schemes \eqref{align:expMPR} and \eqref{align:2stageMagnus} of orders 2 and 4, respectively, with $\beta = 3$ for the initial decay. 
In each time step, we apply the Lanczos process together with the fast algorithm. 
%%%
The left figure gives \eqref{align:timePropagationError} using \eqref{align:expMPR} for fixed $N=2$, $m=7$, and various choices of $K$, revealing the expected order of convergence.  
Additionally, we use asterisks to indicate the perturbation error
\begin{align}\label{align:timePropagationPerturbationError}
 \max_{\bfj\in\Kc} \left|\left(c_{\text{pol}}^n-\bar{c}_{\text{pol}}^n\right)_{\bfj}\right|,
\end{align}
where $\bar{c}_{\text{pol}}^n$ comes from an unperturbed application of Lanczos. 
The latter error becomes eventually dominant, but decreases arbitrarily for increasing $K$. 
%%%
The middle picture shows the case of fixed $N=2$, $K=40$, and various choices of $m$ using \eqref{align:expMPR}. 
Due to the perturbation error being dominant, larger choices of $m$ yield no additional accuracy. 
%%%
Finally, the right picture illustrates the error behavior for fixed $K=40$, $m=7$, and various choices of $N$ using both schemes. 
The order of \eqref{align:2stageMagnus} is not revealed before the perturbation error dominates. 
%%% 
In all cases, \eqref{align:2stageMagnus} has been employed
with $h = 1\mh04$ and 20 unperturbed Lanczos steps in each time step to obtain a reference.
%------------------------
%%%%%%%%%%%%%%%%%%%%%%%%%%%%%%%%%%%%%%%%%%%%%%%%%%%%%%%%%%%%%%%%%%%%%%%%%%%%%%%%%%%%%%
%%%%%%%%%%%%%%%%%%%%%%%%%%%%%%%%%%%%%%%%%%%%%%%%%%%%%%%%%%%%%%%%%%%%%%%%%%%%%%%%%%%%%%
\section*{Conclusion}
We have presented a fast algorithm for the efficient treatment of the coefficient ODE resulting 
from spatial discretization of the linear Schr\"odinger equation in higher dimensions by a spectral Galerkin method. As time discretization of this ODE typically involves products of the time-dependent Galerkin matrix with a vector, 
assembling this matrix and doing the multiplication 
explicitly is prohibitive due to the complexity of the problem -- even with a reduced basis and even more so in each time step. The fast algorithm provides a 
direct approach for this problem to circumvent complexity issues and reduce computational efforts considerably. It consists 
of a sequential, fast application of coordinate matrices formally inserted into the polynomially approximated potential and scales only linearly in the size of the 
chosen basis for any choice of index set reduction. Using a full index set, our quadrature-free procedure is equivalent to 
Gauss--Hermite quadrature with exactly as many nodes as there are basis functions in each direction. 
For a hyperbolically reduced index set, we have analyzed the resulting quadrature and 
index set reduction errors by casting the problem as 
an examination on binary trees. As it turns out, both errors decay rapidly if the underlying potential is sufficiently smoother than the exact solution. 
Possible issues in the context of a Lanczos-based time propagation scheme have been discussed. 
The analysis as given in the present work is limited to the case of the solution being essentially supported within a cube. 
However, we point out that the fast algorithm has been successfully used in combination with a moving wavepacket basis, 
and that it constitutes a generic strategy that can be adopted for spectral discretizations based on 
orthogonal polynomials for linear problems involving boundary conditions other than the Schr\"odinger equation with Hermite functions.\vspace{\baselineskip}\\
%%%%%%%%%%%%%%%%%%%%%%%%%%%%%%%%%%%%%%%%%%%%%%%%%%%%%%%%%%%%%%%%%%%%%%%%%%%%%%%%%%%%%%
%%%%%%%%%%%%%%%%%%%%%%%%%%%%%%%%%%%%%%%%%%%%%%%%%%%%%%%%%%%%%%%%%%%%%%%%%%%%%%%%%%%%%%
\emph{\textbf{Acknowledgement:}} The author is grateful to Christian Lubich for his valuable suggestions and for helpful discussions. Bernd Brumm is funded by the DFG Priority Program 1324 and is associated with DFG Research Training Group 1838.

%%%%%%%%%%%%%%%%%%%%%%%%%%%%%%%%%%%%%%%%%%%%%%%%%%%%%%%%%%%%%%%%%%%%%%%%%%%%%%%%%%%%%%
%%%%%%%%%%%%%%%%%%%%%%%%%%%%%%%%%%%%%%%%%%%%%%%%%%%%%%%%%%%%%%%%%%%%%%%%%%%%%%%%%%%%%%


\begin{thebibliography}{}
\bibitem{AS65} {\scshape M. Abramowitz and I.A.Stegun}, \emph{ Handbook of Mathematical Functions}, Dover, New York, 1972 (10th printing).
\bibitem{AC09} {\scshape G. Avila and T. Carrington}, \emph{ Nonproduct quadrature grids for solving the vibrational Schr\"odinger equation}, J. Chem. Phys., 131 (2009), 174103-1 -- 174103-15.% \verb#doi: 10.1063/1.3246593#.
\bibitem{AC11a} {\scshape \dots}, \emph{ Using nonproduct quadrature grids to solve the vibrational Schr\"odinger equation in 12D}, J. Chem. Phys., 134 (2011), 054126-1 --  054126-16.%\verb#doi: 10.1063/1.3549817#.
\bibitem{AC11b} {\scshape \dots}, \emph{ Using a pruned basis, a non-product quadrature grid, and the exact Watson normal-coordinate kinetic energy operator to solve the vibrational Schrodinger equation for C2H4}, J. Chem. Phys., 135 (2011),  064101-1 -- 064101-12.%\verb#doi: 10.1063/1.3617249#.
\bibitem{AC12} {\scshape \dots}, \emph{ Solving the vibrational Schr\"odinger equation using bases pruned to include strongly coupled functions and compatible quadratures}, J. Chem. Phys., 137 (2012), 174108- 1-- 174108-13.%\verb#doi: 10.1063/1.4764099#.
%\bibitem{AC13} {\scshape \dots}, \emph{ Solving the Schroedinger equation using Smolyak interpolants}, J. Chem. Phys., 139 (2013), \verb#doi: 10.1063/1.4821348#.
%\bibitem{BM00} {\scshape S.~Blanes and P.C.~Moan},  \emph{ Splitting methods for the time-dependent Schr\"odinger equation}, Phys. Lett. A, 265 (2000), pp.~35--42.
\bibitem{BCOR09} {\scshape S. Blanes, F. Casas, J.A. Oteo and J. Ros}, \emph{ The Magnus expansion and some of its applications}, Phys. Rep., 470 (2009), pp.~151--238.
%\bibitem[Blanes \& Moan(2000)]{BM00} {\scshape S.~Blanes and P.C.~Moan}, \emph{ Splitting methods for the time-dependent Schr\"odinger equation}, Physics Letters A, 265 (2000), pp.~35--42.
%\bibitem{BF12} {\scshape D.~Blazevski} \& J.~Franklin, \emph{ Using scattering theory to compute invariant manifolds and numerical results for the laser driven Hénon-Heiles system}, Chaos, 22 (2012), doi: 10.1063/1.4767656.
\bibitem{BG04} {\scshape H.-J. Bungartz and M. Griebel}, \emph{ Sparse Grids}, Acta Numerica, 13 (2004), pp.~147--269.
\bibitem{CQHZ06} {\scshape C. Canuto, A. Quarteroni, M.Y. Hussaini, and T.A. Zang}, \emph{ Spectral Methods: Fundamentals in Single Domains}, Springer, Berlin, 2006.
\bibitem{CJX14} {\scshape Y. Cao, Y. Jiang, and Y. Xu}, \emph{A fast algorithm for orthogonal polynomial expansions on sparse grids}, J. Complexity, 30 (2014), pp.~683--715. %10.1016/j.jco.2014.04.001
%\bibitem{BC93} {\scshape T.~Carrington \& M.J.~Bramley}, \emph{ A General Discrete Variable Method to Calculate Vibrational Energy Levels of Three- and Four-Atom Molecules}, J. Chem. Phys., 99 (1993), pp.~8519--8541.
%\bibitem{CR96} {\scshape T.~Carrington and P.-N.~Roy}, \emph{ A direct-operation Lanczos approach for calculating energy levels}, Chemical Physics Letters, 257 (1996), pp.~98--104.
%\bibitem{C62} {\scshape C.W. Clenshaw}, \emph{ Chebyshev Series for Mathematical Functions}, Mathematical Tables, Vol. 5, H.M. Stationery Office, London, 1962.
\bibitem{FG09} {\scshape E. Faou and V. Gradinaru}, \emph{ Gauss--Hermite wavepacket dynamics: convergence of the spectral and pseudo-spectral approximation}, IMA J. Num. Anal., 29 (2009), pp.~1023--1045.
\bibitem{FGL09} {\scshape E. Faou, V. Gradinaru, and Ch. Lubich}, \emph{ Computing semiclassical quantum dynamics with Hagedorn wavepackets}, SIAM J. Sci. Comp., 31 (2009), pp.~3027--3041.
\bibitem{G11} {\scshape L. Gauckler}, \emph{ Convergence of a split-step Hermite method for the Gross--Pitaevskii equation}, IMA J. Num. Anal., 31 (2011), pp.~396--415.
\bibitem{G12} {\scshape W. Gautschi}, \emph{Numerical Analysis: An Introduction}, Birkh\"{a}user, Basel, 2012 (2nd ed.).
\bibitem{GG98} {\scshape T. Gerstner and M. Griebel}, \emph{ Numerical integration using sparse grids}, Numer. Algor., 18 (1998), pp.~209--232.
\bibitem{G07a} {\scshape V. Gradinaru}, \emph{ Fourier Transform on Sparse Grids: Code Design and Application to the Time Dependent Schr\"odinger Equation}, Computing, 80 (2007), pp.~1--22.
\bibitem{G07b} {\scshape \dots}, \emph{ Strang Splitting for the Time Dependent Schr\"odinger Equation on Sparse Grids}, SIAM J. Num. Anal., 46 (2007), pp.~103--123.
%\bibitem{GH13} {\scshape V.~Gradinaru and G.~Hagedorn}, \emph{ Convergence of a semiclassical wavepacket based time-splitting for the Schr\"odinger equation}, to appear in Numerische Mathematik, \verb#doi:10.1007/s00211-013-0560-6#.%see \verb#www.math.vt.edu/people/hagedorn/gradhag2.pdf#.
\bibitem{HL03} {\scshape M. Hochbruck and Ch.~Lubich}, \emph{ On Magnus integrators for time-dependent Schr\"odinger equations}, SIAM J. Num. Anal., 41 (2003), pp.~945--963.
%\bibitem{IN99} {\scshape A.~Iserles and S.P.~N\o{}rsett}, \emph{ On the solution of linear differential equations in Lie groups}, R. Soc. Lond. Philos. Trans. Ser. A Math. Phys. Eng. Sci., 357 (1999), pp.~983--1019.
%\bibitem{IMNZ00} {\scshape A.~Iserles, H.Z.~Munthe-Kaas, S.P.~N\o{}rsett and A.~Zanna}, \emph{ Lie-group methods}, Acta Numerica, 9 (2000), pp.~215--365.
\bibitem{LC00} {\scshape J.C. Light and T. Carrington}, \emph{ Discrete variable representations and their utilization}, Adv. Chem. Phys., 114 (2000), pp.~263--310. 
\bibitem{L08} {\scshape Ch. Lubich}, \emph{ From Quantum to Classical Molecular Dynamics: Reduced Models and Numerical Analysis}, Europ. Math. Soc., Zurich, 2008.
%\bibitem{K06} {\scshape A.~Klimke}, \emph{ Construction of Hierarchical Polynomial Sparse Grid Interpolants using the Fast Discrete Cosine Transform}, IANS preprint 2006/007, University of Stuttgart, see \verb#http://preprints.ians.uni-stuttgart.de#.
\bibitem{PKM94} {\scshape U. Peskin, R. Kosloff and N. Moiseyev}, \emph{ The Solution of the Time-Dependent Schr\"odinger Equation by the (t, t') Method: The Use of Global Polynomial Propagators for Time-Dependent Hamiltonians}, J. Chem. Phys., 100 (1994), pp.~8849--8855. 
\bibitem{NR07} {\scshape W.H. Press, S.A., Teukolsky, and W.T. Vetterling, B.P. Flannery}, \emph{Numerical Recipes: The Art of Scientific Computing}, Cambridge University Press, New York, 2007 (3rd ed.).
\bibitem{S63} {\scshape S.A. Smolyak}, \emph{ Quadrature and interpolation formulas for tensor products of certain classes of functions}, Dokl. Akad. Nauk SSSR, 4 (1963), pp.~240--243.
%\bibitem{SB02} {\scshape J. Stoer and R. Bulirsch}, \emph{Introduction to Numerical Analysis}, Springer, New York, 2002 (3rd ed.).
\bibitem{Th00} {\scshape B. Thaller}, \emph{ Visual Quantum Mechanics}, Springer, New York, 2000.
\bibitem{VL77} {\scshape Ch.  Van Loan}, \emph{ The Sensitivity of the Matrix Exponential}, SIAM J. Numer. Anal., 14 (1977), pp. 971--981.
\bibitem{Z91} {\scshape Ch. Zenger}, \emph{Sparse Grids}, in: Parallel algorithms for partial differential equations, Notes on Numerical Fluid Mechanics 31, W. Hackbusch, ed., Vieweg, Braunschweig, 1991, pp.~241--251.
\end{thebibliography}
\end{document}